\documentclass[11pt]{amsart}

\usepackage[T1]{fontenc}
\usepackage{lmodern}

\usepackage[left=2.5cm,right=2cm,top=2cm,bottom=2cm]{geometry}

\usepackage{amssymb, amsmath, amsthm}
\usepackage{mathabx}
\usepackage{graphicx}
\usepackage[utf8]{inputenc}
\usepackage{hyperref}
\usepackage{subcaption}
\usepackage{easyReview}

\usepackage{young}
\usepackage{mathtools}
\usepackage{rotating}
\usepackage{tabularx}

\usepackage{enumitem}

\newtheorem{thm}{Theorem}[section]
\newtheorem*{thm*}{Theorem}

\newtheorem{lem}[thm]{Lemma}
\newtheorem{prop}[thm]{Proposition}
\theoremstyle{definition}

\newtheorem{rem}[thm]{Remark}
\newtheorem{ex}{Example}
\newcommand{\R}{\mathbb{R}}
\newcommand{\Z}{\mathbb{Z}}
\newcommand{\N}{\mathbb{N}}
\newcommand{\tq}{\colon}
\newcommand{\weyl}{\mathcal{W}}
\newcommand{\schub}{\mathcal{S}}
\newcommand{\wmax}{\widetilde w}
\newcommand{\wprimemax}{\widetilde w'}
\newcommand{\type}{\mathrm{type}}
\newcommand{\invp}{\mathrm{Inv}^{+}}
\newcommand{\invn}{\mathrm{Inv}^{-}}
\renewcommand{\mod}{\ \mathrm{mod}\ }
\newcolumntype{C}[1]{>{\centering\arraybackslash}p{#1}}
\newcommand{\isB}{\mathsf{isB}}
\newcommand{\typew}{\mathsf{Type}}

\allowdisplaybreaks %Permite dividir equação em páginas diferentes

\begin{document}

\title{Integral homology of real isotropic and odd orthogonal Grassmannians}

%% arXiv

\author{Jordan Lambert}
\address[Jordan Lambert]{Department of Mathematics -- ICEx, Universidade Federal Fluminense, Volta Redonda 27213-145, Rio de Janeiro, Brazil}
\email[Corresponding author]{jordanlambert@id.uff.br}

\author{Lonardo Rabelo}
\address[Lonardo Rabelo]{Department of Mathematics, Federal University of Juiz de Fora, Juiz de Fora 36036-900, Minas Gerais, Brazil}
\email{lonardo@ice.ufjf.br}

\thanks{Supported by FAPESP grant number 13/10467-3 and 14/27042-8, and the Coordination for the Improvement of Higher Level Personnel -- Capes.}

% \subjclass is required.

\keywords{Grassmannian permutation; Isotropic Grassmannian; Integral homology}

\subjclass[2010]{Primary 05A05, 05E15, 14M15, 57T15}

\begin{abstract}
We obtain a combinatorial expression for the coefficients of the boundary map of real isotropic and odd orthogonal Grassmannians providing a natural generalization of the formulas already obtained for Lagrangian and maximal isotropic Grassmannians. The results are given in terms of the classification into four types of covering pairs among the Schubert cells when identified with signed $k$-Grassmannian permutations. It turns out that these coefficients only depend on the positions changed over each pair of permutations. As an application, we give an orientability criterion, exhibit a symmetry of these coefficients and, compute low-dimensional homology groups.
\end{abstract}

\maketitle

\section{Introduction}

In this paper, we provide an explicit formula for the coefficients of the boundary map of odd orthogonal and isotropic real Grassmannians. This work corresponds to a natural generalization of the results of the first author paper \cite{Rab16} where these coefficients were given in the context of the Lagrangian and maximal isotropic Grassmannians. It provides a contribution to the study of the topology of real flag manifolds which is a very delicate object in comparison with the complex ones since the latter does not have torsion in homology. The type of formulas we have found is compatible with those of Burghelea--Hangan--Moscovici--Verona \cite{BHMV73} which gave the first results about the topology of real flag manifolds in the 1970s. 

We consider the Grassmannian manifolds as minimal flag manifolds of $G$ given either by the odd orthogonal group $\mathrm{SO}(n,n+1)$ or by the symplectic group $\mathrm{Sp}(n,\R)$ with Lie algebra $\mathfrak{g}$ given respectively by $\mathfrak{so}(n,n+1)$ or $\mathfrak{sp}(n,\R)$. We denote the root system by $\Sigma=\{a_0,\ a_1=\varepsilon_{2}-\varepsilon_{1}\ , \dots,\ a_{n-1}=\varepsilon_{n}-\varepsilon_{n-1}\}$ such that $a_{0}=\varepsilon_{1}$ for type $B$ and $a_{0}=2\varepsilon_{1}$ for type $C$. The  Grassmannians are the minimal flag manifolds $G/P_{(k)}$, $0\leq k \leq n-1$, where $P_{(k)}$ is the parabolic subgroup corresponding to the maximal proper subset $(k)=\Sigma-\{a_{k}\}$. For each choice of group $G$, it corresponds to the following  geometric realization: for $G= \mathrm{SO}(n,n+1)$, we have the orthogonal Grassmannian $\mathrm{OG}(n-k, 2n+1)$, the set of $(n-k)$-dimensional isotropic subspaces in the vector space $V=\R^{2n+1}$ equipped with a nondegenerate symmetric bilinear form; and, for $G=\mathrm{Sp}(n,\R)$, we have the isotropic Grassmannian $\mathrm{IG}(n-k,2n)$, the set $(n-k)$-dimensional isotropic subspaces in the symplectic vector space $V=\R^{2n}$. 

The cellular structure of these Grassmannians will be given by the Schubert cells parametrized by the minimal representatives $\weyl_n^{(k)}$ of Weyl group modulo the subgroup generated by the reflections of $(k)$. We use a permutation model to identify cells with the set of signed $k$-Grassmannians permutations, i.e., signed permutations of the form $w = u_1\, \cdots\, u_k | \overline{\lambda_r}\, \cdots\, \overline{\lambda_1}\, v_{1}\, \cdots \, v_{n-k-r}$, where $0\leqslant r\leqslant n-k$, $0<u_1<\cdots <u_k, 0<\lambda_1 < \cdots < \lambda_r$ and $0<v_1 <\cdots < v_{n-k-r}$. This approach has been successfully employed over related problems in Schubert calculus of such Grassmannians (see Pragacz-Ratajski \cite{PR96}, Buch-Kresch-Tamvakis \cite{BKT17} and Ikeda-Matsumura \cite{Ike15}). Our setting is based on the paper of the same authors \cite{LR19} where the coverings pairs $w,w'\in \weyl_{n}^{(k)}$ -- those satisfying $w'\leq w$ in the Bruhat-Chevalley order with $\ell(w)=\ell(w')+1$ -- are classified into four types (B1, B2, B3 and B4). We associate with each permutation $w \in \weyl_{n}^{(k)}$ a double partition $(\alpha,\lambda)$ for which we have the corresponding half-shifted Young diagram (HSYD).
Since one can associate with each covering pair $w,w'$ the corresponding double partitions $(\alpha,\lambda)$ and $(\alpha',\lambda')$ (see Proposition \ref{prop:doublepartitiontypes}), our main results have a very nice representation in terms of these diagrams facilitating its comprehension.

According to ---,San Martin \cite{RS19}, all data required to determine the coefficients of the homology groups is extracted from the set $\Pi_w$ of positive roots sent to negative by $w^{-1}$. The knowledge of the set $\Pi_{w}$ is of high importance since its cardinality gives the dimension of the Schubert varieties (for split real forms). A great source of inspiration was the papers of Ikeda-Naruse \cite{IN08} and Graham-Kreiman \cite{GK15} where $\Pi_w$ is described with aid of the shifted Young diagrams in the context of the maximal isotropic Grassmannians. We show here how to generalize these ideas for all Grassmannians of types B and C by using the HSYD's associated to the double partition $(\alpha,\lambda)$ of $w$. An essential idea consists in the determination of $\Pi_w$ as a set of inversions of the corresponding permutation such that each box of the diagram contains a uniquely defined inversion of $\Pi_w$ (cf. Sections \ref{subsec:doublepart} and \ref{subsec:inversions}).

Once the set of inversions $\Pi_w$ is known and the covering pairs $w,w'$ are classified, by (\cite{RS19}, Theorem 2.8), it remains to compute the difference $\sigma(w) - \sigma(w')$, where $\sigma(w)$ is the sum over the roots of $\Pi_w$. Since this is a multiple $\kappa(w,w')$ of a unique root by Equation \eqref{eq:kappa}, the coefficient $c(w,w')$ in the boundary map will be either $0$ or $\pm 2$ depending on the parity of $\kappa(w,w')$. Our most outstanding result given by Theorem \ref{thm:mainthm} is that $c(w,w')$ is completely determined by a well-defined triple $(P,T,Q)$ composed of the positions changed in permutations $w$ and $w'$. Besides its simplicity, it represents an absolutely new interpretation of these coefficients. This is the content of Theorem \ref{thm:dualheight} which relates the computation of $\kappa$ with the height of a given root. 

As a direct consequence of our formulas, a criterion of orientability is easily obtained (see Proposition \ref{prop:orient}). Furthermore, we show that the incidence graph whose vertices are the permutations of $\weyl_{n}^{(k)}$ and the edges ``$\rightarrow$'' and ``$\Rightarrow$'' corresponds to $0$ and $\pm 2$ coefficients, respectively, admits an elegant symmetry (see Figure \ref{fig:n4k2_incidence} and Proposition \ref{prop:dualtype}). We also finish by computing some low dimensional homology groups (see Section \ref{subsec:lowdim}). These provide a combinatorial description of the results in this specific context related to the more general ones as obtained by Patr\~ao-San Martin-Santos-Seco \cite{PS12} concerning the orientability and del Barco-San Martin \cite{BS19} with respect to the low dimensional groups.

The article is organized as follows.
Section \ref{sec:grassm} describes the isotropic and odd orthogonal Grassmannian, some properties of covering relations, and the half-shifted Young diagrams.
Section \ref{sec:boundarymap} describes the boundary maps of the cellular homology of Grassmannians of type B and C,  the set of inversions applied to such Grassmannians, and it states the main result in this work. Some additional results about orientability and symmetry are also proven. Section \ref{sec:proofs} is devoted to prove the main result. We conclude with Section \ref{sec:finalcomments} where we provide some perspectives for future work.

\section{Isotropic and odd orthogonal Grassmannians}\label{sec:grassm}

We let $\N=\{1,2,3, \dots\}$ and $\Z$ be the set of integers. For $n,m\in\Z$, where $n\leqslant m$, denote the set $[n,m]=\{n,n+1, \dots, m\}$. For $n\in \N$, denote $[n]=[1,n]$.

Let $G$ be a non-compact semi-simple Lie group. Denote by $\Pi$ a set of roots related to the Lie algebra $\mathfrak{g}$, $\Pi^{\pm}$ the set of positive and negative roots. Fix a simple root system $\Sigma\subset \Pi$. Let $P$ be the minimal parabolic subgroup with Lie algebra $\mathfrak{p}$. We call $\mathbb{F} = G/P$ the maximal flag manifold of $G$. The Weyl group $\weyl$ of $G$ is the product of simple reflections $s_{i}=s_{a_{i}}$ through simple roots $a_{i}\in\Sigma$. The length $\ell(w)$ of $w \in \mathcal{W}$ is the number of simple reflections in any reduced decomposition of $w$.  There is a partial order in the Weyl group called the Bruhat-Chevalley order: we say that $w_{1}\leqslant w_{2}$ if given a reduced decomposition $w_2 = s_{j_{1}} \cdots s_{j_{r}}$ then $w_{1}=s_{j_{i_{1}}}\cdots s_{j_{i_{k}}}$ for some $1\leqslant i_1\leqslant \cdots \leqslant i_r\leqslant r$. It is known that $\mathcal{W}$ has a maximum element $w_0$ which is an involution, i.e, $w_0^2=1$.

A subset of simple roots $\Theta\subset \Sigma$ is associated with the parabolic subgroup $P_\Theta$ in $G$, which contains $P$. The corresponding flag manifold $\mathbb{F}_\Theta = G/P_{\Theta}$ is called a partial flag manifold of $G$. The subgroup $\mathcal{W}_\Theta$ is the subgroup of the Weyl group $\weyl$ generated by the reflections with respect to the roots $\alpha \in \Theta$.
We also define the subset $\mathcal{W}^{\Theta}$ of $\mathcal{W}$ by $\mathcal{W}^{\Theta} = \{w \in \mathcal{W} \colon \ell(ws_a) = \ell(w) + 1 ~,~\forall a \in \Theta \}$. Since there exists a unique element $w^{\Theta} \in \mathcal{W}^{\Theta}$ of minimal length in each coset $w\mathcal{W}_{\Theta}$,  $\mathcal{W}^{\Theta}$ is called the subset of minimal representatives of the cosets of $\mathcal{W}_{\Theta}$ in $\mathcal{W}$. 

The Bruhat decomposition presents the flag manifolds as disjoint union $\mathbb{F}_{\Theta }=\coprod_{w\in \mathcal{W}^{\Theta }} \ N\cdot wb_{\Theta}$ where $N$ is the nilpotent subgroup in the Iwasawa decomposition of $G$. A Schubert variety is the closure of a Bruhat cell, i.e.,  $\mathcal{S}_{w}=\mathrm{cl}(N\cdot wb_{\Theta })$. The choice of a minimal representative $w\in \weyl^{\Theta}$ gives $\dim(\schub_{w}) = \ell(w)$ since we are in a split case. The Bruhat-Chevalley order defines an order between the Schubert varieties by $\mathcal{S}_{w_{1}}\subset \mathcal{S}_{w_{2}} \mbox{ if, and only if, } w_{1}\leqslant w_{2}$.

First, consider the symplectic group $\mathrm{Sp}(n,\mathbb{R})$ with Lie algebra $\mathfrak{sp}(n,\mathbb{R})$ of type $C$. The root system of type $C$ is realized as a set of vectors $\Pi=\{\pm \varepsilon_{i}\pm \varepsilon_{j} \tq 1\leqslant i< j\leqslant n\}\cup \{\pm 2\varepsilon_{i} \tq 1\leqslant i\leqslant n\}$ in the Euclidean space $\mathbb{R}^{n}=\oplus_{i=1}^{n}\mathbb{R}\varepsilon_{i}$. Denote the (positive) simple roots by $a_{0}=2\varepsilon_{1}$ and $a_{i}=\varepsilon_{i+1}-\varepsilon_{i}$ for $1\leqslant i < n$. Then, the set of all positive roots is $\Pi^{+}=\{\varepsilon_{j} \pm \varepsilon_{i} \tq 1\leqslant i< j\leqslant n\}\cup \{2\varepsilon_{i} \tq 1\leqslant i \leqslant n\}$. Given $k\in [0,n-1]$, the minimal flag manifold $\mathrm{IG}(n-k,2n)=\mathrm{Sp}(n,\mathbb{R})/P_{(k)}$, for $(k)=\Sigma-\{a_{k}\}$, is called a \emph{isotropic Grassmannian} since it parametrizes $(n-k)$-dimensional isotropic subspaces of a real $2n$-dimensional symplectic vector space.

Now, consider the orthogonal group $\mathrm{SO}(n,n+1)$ with Lie algebra $\mathfrak{so}(n,n+1)$ of type $B$. The root system of type $B$ is realized as a set of vectors $\Pi=\{\pm \varepsilon_{i}\pm \varepsilon_{j} \tq 1\leqslant i< j\leqslant n\}\cup \{\pm \varepsilon_{i} \tq 1\leqslant i\leqslant n\}$ in the Euclidean space $\mathbb{R}^{n}=\oplus_{i=1}^{n}\mathbb{R}\varepsilon_{i}$. Denote the (positive) simple roots by $a_{0}=\varepsilon_{1}$ and $a_{i}=\varepsilon_{i+1}-\varepsilon_{i}$ for $1\leqslant i < n$. Then, the set of all positive roots is $\Pi^{+}=\{\varepsilon_{j} \pm \varepsilon_{i} \tq 1\leqslant i< j\leqslant n\}\cup \{\varepsilon_{i} \tq 1\leqslant i \leqslant n\}$. Given $k\in [0,n-1]$, the minimal flag manifold $\mathrm{OG}(n-k,2n+1)=\mathrm{SO}(n,n+1)/P_{(k)}$, for $(k)=\Sigma-\{a_{k}\}$, is called a \emph{odd orthogonal Grassmannian} since it parametrizes $(n-k)$-dimensional isotropic subspaces of a real $(2n+1)$-dimensional vector space equipped with a nondegenerate symmetric bilinear form.

\subsection{Permutation model}\label{subsec:model}

For both types of Grassmannians, we denote $s_{i}=s_{a_{i}}$, for $i \in [0,n-1]$, the simple reflection given by the simple root $a_{i}$. The Weyl group $\weyl_{n}$ for root systems $B_n$ and $C_n$, also called \emph{hyperoctahedral group}, is the semidirect product $S_{n}\ltimes \Z^{n}_{2}$. We also realize as the set of permutations of $S_n$ with a sign (plus or minus) attached to each entry; we will write these elements as barred permutations using the bar to denote a negative sign, and we take the natural order on them. Then a permutation $w \in \weyl_{n}$, usually denoted in one-line notation $w=w(1)\ w(2)\cdots w(n)$, satisfies the relation $\overline{w(i)}=w(\overline{i})$. With respect to this realization, the length of $w \in \mathcal{W}_{n}$ is given by the following formula (\cite{BB05},Eq.(8.3))
\begin{equation}\label{eq:length}
\ell(w) = \mathrm{inv}(w(1),\ldots, w(n)) - \sum_{ \mathclap{\{j\mid w(j)<0\}}}\ w(j)
\end{equation}
where 
\begin{equation*}
\mathrm{inv}(w(1),\ldots,w(n)) = \# \{(i,j) \colon 1\leq i<j\leq n, w(i)>w(j)\}. 
\end{equation*}

The simple reflections $s_{0},\dots, s_{n-1}$ act on the right of a permutation $w$ in $\weyl_{n}$ by
\begin{align*}
w(1)\ w(2) \cdots w(n) \cdot s_{0} &=\overline{w(1)}\ w(2) \cdots w(n); \\
w(1)\cdots w(i)w(i+1) \cdots w(n)\cdot s_{i} &=w(1)\cdots w(i+1) w(i) \cdots w(n) \mbox{ , } 1\leqslant i<n
\end{align*}

The hyperoctahedral group $\weyl_{n}$ is also the Coxeter group of type B generated by $s_i$ and subject to the relations: (i) $s_i^2=1$, for $i\geq 0$; (ii) $s_0s_1s_0s_1 = s_1s_0s_1s_0$; (iii) $s_{i+1}s_i s_{i+1} = s_i s_{i+1}s_i$, for $i\in [n-1]$; (iv) $s_is_j = s_js_i$, for $|i-j|\geqslant 2$.

As defined above, for $(k)=\Sigma-\{a_{k}\}$, the corresponding subgroup $\mathcal{W}_{(k)}$ is generated by $s_i$, with $i\neq k$. Notice that $\mathcal{W}_{(k)}\cong \mathcal{W}_k \times S_{n-k}$, where $\mathcal{W}_k$ is the subgroup generated by $s_i$, $i\in[0,k]$. Define by $\smash{\mathcal{W}^{(k)}_{n}}\subset \mathcal{W}_n$ the set of minimal length coset representatives of $\weyl_{n}/\mathcal{W}_{(k)}$, which parametrizes the Schubert varieties in $\mathrm{IG}(n-k,2n)$ and $\mathrm{OG}(n-k,2n+1)$. This indexing set $\mathcal{W}^{(k)}_{n}$ can be identified by a set of signed permutations of the form
\begin{equation}\label{eq:kgrass}
w = w_{u,\lambda} = u_1\, \cdots\, u_k | \overline{\lambda_r}\, \cdots\, \overline{\lambda_1}\, v_{1}\, \cdots \, v_{n-k-r}
\end{equation}
where $0\leqslant r\leqslant n-k$ and
\begin{align}\label{eq:kgrass_conds}
0<u_1<\cdots <u_k\ , &&& u_i=w(i), \mbox{ for } i \in [k]; \nonumber\\
0<\lambda_1 < \cdots < \lambda_r\ ,  &&& \overline{\lambda_i} = w(k+r-i+1), \mbox{ for } i \in [r]; \\
0<v_1 <\cdots < v_{n-k-r}\ ,  &&& v_i = w(k+r+i), \mbox{ for } i\in [n-k-r]. \nonumber
\end{align}

They are called signed $k$-Grassmannian permutations. 
The longest element $w_{0}^{k} \in \weyl_{n}^{(k)}$ is the $k$-Grassmannian permutation given by 
\begin{equation}\label{eq:longestw}
w_{0}^{k}=1\,2\,\cdots\, k| \overline{n\vphantom{1}}\,\overline{n-1}\, \cdots\, \overline{k+1}.
\end{equation}

Throughout the text, we assign $r$ as the length of $\lambda$, i.e., $\ell(\lambda)=r$.

\subsection{Bruhat order and covering relations}\label{subsec:cover}
We recall some results of \cite{LR19} about the Bruhat-Chevalley order of $\weyl_{n}^{(k)}$ in the permutation model. Let $w,w'  \in \mathcal{W}_{n}$ with $w'\leqslant w$ and $\ell(w)=\ell(w')+1$, i.e., if $w=s_{i_{1}}\!\cdots s_{i_{\ell(w)}}$ is a reduced decomposition then $w'=s_{i_{1}}\!\cdots \widehat{s_{i_{j}}}\! \cdots s_{i_{\ell(w)}}$ is a reduced decomposition as well. In this case, we say that it is a covering relation where $w$ covers $w'$. According to \cite{LR19}, there is a classification of pairs $w,w'$ where $w$ covers $w'$ in $\weyl_{n}^{(k)}$. We now define the types of existing pairs. Suppose that $w$ and $w'$ are permutations in $\weyl_{n}^{(k)}$ according to Equation \eqref{eq:kgrass_conds}. 

%We say that \add{two permutations $w,w'\in \weyl_{n}^{(k)}$ are of type} \remove{$\typew(w,w')$ is} B1,B2, B3 or B4 accordingly to the following description for $w$ and $w'$:

We say that $w,w'$ is \emph{a pair of type B1} if they are written as follows:
\begin{equation*}
w = \cdots\, |\, \cdots\,  \overline{1}\,  \cdots \,  \quad\mbox{ and} \quad
w' = \cdots\, |\, \cdots\,  1\,  \cdots.
\end{equation*}
In other words, if $w$ is such that $\lambda_{1}=1$, then $w'$ is obtained from $w$ by removing the negative sign from $\overline{1}$. 

We say that $w,w'$ is \emph{a pair of type B2} if they are written as follows:
\begin{equation*}
w = \cdots\, |\, \cdots\,  \overline{a}\,  \cdots \, (a-1)\, \cdots \, \quad \mbox{ and} \quad
w' = \cdots\, |\, \cdots\, \overline{a-1}\,  \cdots \, a\, \cdots,
\end{equation*}
where $a>0$. In other words, there are $t \in [r]$ and $q \in [n-k-r]$ such that $\lambda_{t}=a$ and $v_{q}=a-1$, and $w'$ is obtained from $w$ by switching $v_{q}$ and $\lambda_{t}$.

We say that $w,w'$ is \emph{a pair of type B3} if they are written as follows:
\begin{equation*}
w = \cdots\, a \,\cdots\, |\,  \cdots \, (a-x)\, \cdots \, \quad \mbox{ and} \quad
w' = \cdots\,  (a-x) \,\cdots \,|\, \cdots\, a\, \cdots,
\end{equation*}
where $a > x > 0$. In other words, there are $p \in [k]$ and $q \in [n-k-r]$ such that $u_{p}=a$ and $v_{q}=a-x$. The permutation $w'$ is obtained from $w$ by switching $u_{p}$ and $v_{q}$.

Finally, we say that $w,w'$ is \emph{a pair of type B4} if they are written as follows:
\begin{equation*}
w = \cdots\, (a-x) \,\cdots\, |\,  \cdots \, \overline{a}\, \cdots \, \quad \mbox{ and} \quad
w' = \cdots\,  a \,\cdots \,|\, \cdots\, \overline{a-x}\, \cdots,
\end{equation*}
where $a > x  > 0$. In other words, there are $p \in [k]$ and $t \in [r]$ such that $u_{p}=a-x$ and $\lambda_{t}=a$. The permutation $w'$ is obtained from $w$ by switching $u_{p}$ and $\lambda_{t}$. 

We will denote the type of a pair by $\typew(w,w')$. For each type of pair $w,w'$, also consider the integers $P=P(w,w')$, $T=T(w,w')$, and $Q=Q(w,w')$ which correspond to the positions where $w$ changes when compared to $w'$. They can be represented as a right action in the complete notation as following:
\begin{equation}
\begin{array}{cl}\label{eq:rightaction}
\mbox{(B1)} & w = w'\cdot (T,\overline{T})\\
\mbox{(B2)} & w = w'\cdot (\overline{T},Q)(\overline{Q},T)\\
\mbox{(B3)} & w = w'\cdot (P,Q)(\overline{P},\overline{Q})\\
\mbox{(B4)} & w = w'\cdot (\overline{P},T)(\overline{T},P)
\end{array}
\end{equation}

Precisely, the integers $P$, $T$, and $Q$ can be given by as in Table \ref{tbl:ptq}.
\begin{table}[ht]
\caption{Integers $P$, $T$, and $Q$ for a pair $w,w'$}
\label{tbl:ptq}
\centering
\def\arraystretch{1.2}
\small
\begin{tabular}{|c|c|c|c|}
\hline 
$\typew(w,w')$ & $P$ & $T$ & $Q$ \\ 
\hline \hline
B1 &  & $w^{-1}(\overline{1})=k+r$ &  \\ 
\hline 
B2 &  & $w^{-1}(\overline{\lambda_{t}})=k+r-t+1$ & $w^{-1}(v_{q})=k+r+q$ \\ 
\hline 
B3 & $w^{-1}(u_{p})=p$ &  & $w^{-1}(v_{q})=k+r+q$ \\ 
\hline 
B4 & $w^{-1}(u_{p})=p$ & $w^{-1}(\overline{\lambda_{t}})=k+r-t+1$ &  \\ 
\hline 
\end{tabular}
\end{table}

There is a straight relationship between types of pairs and covering relations.

\begin{thm}[\cite{LR19}, Thm. 5]\label{thm:covering}
Let $w,w'\in \weyl_{n}^{(k)}$. Then $w$ covers $ w'$ if, and only if, $\typew(w,w')$ is B1, B2, B3 or B4.
\end{thm}

\begin{rem} In  Ikeda \cite{Ike15} (Lemmas 3.1 and 3.2) there is a description of the covering relation with respect to the weak Bruhat order (cf. Lemmas 3.1 and 3.2). The above result refines it to the strong Bruhat order.
%\remove{This refines the covering pairs as described in} \cite{Ike15} \remove{(Lemmas 3.1 and 3.2) where it's given with respect to the weak Bruhat order}.
For another approach of the covering relation in terms of $k$-strict partitions we refer to Tamvakis-Wilson \cite{TW16}. 
\end{rem}

\begin{ex}\label{explo:1}[\cite{LR19}, Example 3.4] 
Consider $w=2\,6|\overline{7}\,\overline{5}\,\overline{1}\,3\,4$ where $n=7$ and $k=2$. Theorem \ref{thm:covering} guarantees that there are five permutations covered by $w$:
\begin{align*}
\mbox{ (B1) } & w=2\,6|\overline{7}\,\overline{5}\,\mathbf{\overline{1}}\,3\,4 \,,\, w'_{1}=2\,6|\overline{7}\,\overline{5}\,\mathbf{1}\,3\,4 \mbox{ and } T=5;\\
\mbox{ (B2) } & w=2\,6|\overline{7}\,\mathbf{\overline{5}}\,\overline{1}\,3\,\mathbf{4} \,,\,w'_{2}= 2\,6|\overline{7}\,\mathbf{\overline{4}}\,\overline{1}\,3\,\mathbf{5} \mbox{ and } (T,Q)=(4,7); \\
\mbox{ (B3) } & w=2\,\mathbf{6}|\overline{7}\,\overline{5}\,\overline{1}\,3\,\mathbf{4} \,,\, w'_{3}=2\,\mathbf{4}|\overline{7}\,\overline{5}\,\overline{1}\,3\,\mathbf{6} \mbox{ and } (P,Q)=(2,7) ; \\
\mbox{ (B4) } &w=2\,\mathbf{6}|\mathbf{\overline{7}}\,\overline{5}\,\overline{1}\,3\,4 \,,\, w'_{4}=2\,\mathbf{7}|\mathbf{\overline{6}}\,\overline{5}\,\overline{1}\,3\,4 \mbox{ and } (P,T)=(2,3); \\
\mbox{ (B4) } &w=\mathbf{2}\,6|\overline{7}\,\mathbf{\overline{5}}\,\overline{1}\,3\,4 \,,\, w'_{5}= \mathbf{5}\,6|\overline{7}\,\mathbf{\overline{2}}\,\overline{1}\,3\,4 \mbox{ and } (P,T)=(1,4) .
\end{align*}
\end{ex}

\subsection{Double partitions and half-shifted Young diagrams}\label{subsec:doublepart}

In this section, we will understand how $k$-Grassmannian permutation gives rise to a double partition. In this way, there is a bijection between the set $\mathcal{W}_{n}^{(k)}$ and the set of all double partitions.  Moreover, we will show how each double partition may be represented by a half-shifted Young diagram.

Given $n,k$ integers such that $0\leqslant k < n$, we say that $\Lambda=(\alpha,\lambda)$ is a double partition when $\alpha=(0\leqslant \alpha_{1}\leq \cdots \leq \alpha_{k}\leqslant n-k)$ is a partition and $\lambda=(0<\lambda_{1}<\cdots<\lambda_{r} \leqslant n)$ is a strict partition (if $k=0$ then $\alpha$ is represented as an empty set; an empty $\lambda$ is represented by $r=0$).

Firstly, let us define the diagrams. We define the Young diagram associated to the partition $\alpha$ by
\begin{equation}\label{eq:ydiagram}
\mathcal{D}_{\alpha}=\{(i,j) \in \mathbb{Z}^2\colon 1\leq i\leq k \, , \, 1\leq j\leq \alpha_{i}\}.
\end{equation}

Furthermore, we define the shifted Young diagram associated to the strict partition $\lambda$ by
\begin{equation}\label{eq:sydiagram}
\mathcal{SD}_{\lambda}=\{ (i,j) \in \mathbb{Z}^2\colon 1\leq i \leq r\ \, , \, i\leq j \leq i-1+\lambda_{r+1-i} \}.
\end{equation}

We can insert the diagram $\mathcal{D}_{\alpha}$ of a partition $\alpha$ into a rectangle of dimensions $k\times (n-k)$ while the shifted Young diagram $\mathcal{SD}_{\lambda}$ of a strict partition $\lambda$ fits inside a stair shaped triangle with $n$ lines. Let us denote by $\mathcal{D}_{k,n}$ the set of all partitions whose respective Young diagrams are inside a rectangle $k\times (n-k)$ and by $\mathcal{SD}_{n}$ the set of all strict partitions whose respective shifted Young diagrams are inside a stair shaped diagram of length $n$. We denote by $\mathcal{P}(k,n)$ the set of the pairs $(\alpha,\lambda)$ with $\alpha \in \mathcal{D}_{k,n}$ and $\lambda \in \mathcal{SD}_n$ with the property $\ell(\lambda)\leq \alpha_1$, i.e., 
\begin{align}
0 \leq \alpha_{1} \leq  \cdots \leq \alpha_{k} & \leq  n-k \, ;\notag \\
0 < \lambda_{1} < \cdots < \lambda_{r} & \leqslant n \, ;\label{eq:doublepartition} \\
\ell(\lambda) & \leq  \alpha_{1}. \notag
\end{align} 

Notice in this definition that the strict partition $\lambda$ is empty when $\alpha_{1}=0$.

\begin{figure}[ht]
\centering
\includegraphics[scale=0.7]{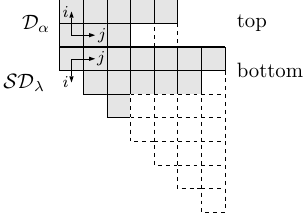}
\caption{A model of a HSYD obtained as a juxtaposition of the diagrams $\mathcal{D}_{\alpha}$ and $\mathcal{SD}_{\lambda}$.}
\label{fig:fig0}
\end{figure}

A half-shifted Young Diagram (HSYD) of the double pair $\Lambda=(\alpha,\lambda) \in \mathcal{P}(k,n)$ is obtained by the juxtaposition of the diagrams $\mathcal{D}_{\alpha}$ and $\mathcal{SD}_{\lambda}$ such that $\mathcal{D}_{\alpha}$ is above $\mathcal{SD}_{\lambda}$. We say that $\mathcal{D}_{\alpha}$ is the top diagram and that $\mathcal{SD}_{\lambda}$ is the bottom diagram. The condition $\alpha_{1}\geq \ell(\lambda)$ is equivalent to say that the number of lines of the bottom diagram is at most the number of full length columns in $\alpha$. Figure \ref{fig:fig0} shows the HSYD for the pair $\alpha = (3,5)$ and $\Lambda = (1,5,7)$.

Now, given a $k$-Grassmannian permutation in one-line notation as in Equation \eqref{eq:kgrass}, let us define its corresponding pair of double partitions $\alpha$ and $\lambda$. The strict partition is given by the negative part of $w$, i.e., $\lambda =(\lambda_{r}>\cdots >\lambda_{1}>0)$.  The partition $\alpha$ is defined by $\alpha_i = u_i - i + d_i$, for $i \in [k]$, where $d_i = \# \{\lambda_{j} \mid \lambda_j > u_i\}$. As shown in \cite{LR19}, Equation (2.10), we have
\begin{equation}\label{eq:alpha3}
\alpha_i = n-k - \mu_{i},
\end{equation}
where $\mu_{i}=\mu_{i}(w)=\#\{v_{j}\mid v_j>u_i\}$, for $i\in [k]$. Since $\mu_{i}$ can be written as $\mu_{i}=n-k-r-\#\{v_{j}\mid v_j<u_i\}$, then
\begin{equation}\label{eq:alpha2}
\alpha_{i}= r+\#\{v_{j}\mid v_j<u_i\}
\end{equation}
and $\alpha$ satisfies $n-k\geq \alpha_k\geq\alpha_{k-1}\geq \cdots \geq \alpha_1\geq0$. 

The partition $\alpha$ counts the number of inversions of the first $k$ entries while the strict partition $\lambda$ the number of remaining inversions given by the negative entries.  Hence, given $w \in \mathcal{W}_{n}^{(k)}$, the length $\ell(w)$ of $w$ is the sum of entries of the pair $(\alpha,\lambda)$, i.e., 
$\ell(w) = |\alpha|+|\lambda|$.

The $k$-Grassmannian permutations are parametrized by the HSYD's.

\begin{lem}[\cite{PR96}, Lemma 1.2] There is a bijection between $\weyl_{n}^{(k)}$ and $\mathcal{P}(k,n)$.
\end{lem}

Given a partition $\alpha=(0\leqslant \alpha_{1}\leq \cdots \leq \alpha_{k}\leqslant n-k)$, we have the \emph{conjugate partition} $\alpha^{*}$ of $\alpha$ which is also a partition defined by $\alpha_{i}^{*}=\#\{\alpha_{j} \tq \alpha_{j}\geqslant i\}$, for all $i \in [n-k]$ satisfying $k\geqslant \alpha^{*}_{1}\geq \cdots \geq \alpha^{*}_{n-k}\geqslant 0$. When we consider the Young diagram $\mathcal{D}_{\alpha}$ of $\alpha$, the conjugate $\alpha^{*}$ is the number of boxes in each column.

If we denote $\mu_{i}^{*} = \#\{u_{j} \tq u_{j}>v_{i}\}$ for $i \in [n-k-r]$, then we have an explicit formula for $\alpha^{*}$ in terms of $\mu_{i}^{*}$. 

\begin{lem}\label{lema:conjugationExplict}
The conjugate partition is given by
\begin{equation}
\alpha_{i}^{*}=
\left\{
\begin{array}{cl}
k & \mbox{if } 1\leqslant i\leqslant r; \\
\mu_{i-r}^{*} & \mbox{if } r+1 \leqslant i \leqslant n-k.
\end{array} 
\right.
\end{equation}
\end{lem}
\begin{proof}
Suppose firstly that $1\leqslant i\leqslant r$. We clearly have $\alpha_{i}^{*}=k$ since $\alpha_{j}\geqslant r$ for every $j$. 
Now, suppose that $r+1 \leq i \leq n-k$. By Equation \eqref{eq:alpha2}, $\#\{v_{l}\tq v_l<u_j\}= \alpha_{j}-r$ for $1\leqslant j \leqslant k$, which is equivalent to $v_1<v_{2}<\cdots < v_{\alpha_{j}-r}<u_j$. Hence, $\alpha_{j}\geqslant i$ if, and only if, $v_{i-r}<u_{j}$. Hence, $\alpha_{i}^{*}= \# \{\alpha_{j} \tq \alpha_j\geq i \, \} = \#\{ u_{j} \tq v_{i-r} < u_j \} = \mu_{i-r}^{*}$.
\end{proof}

We now present an adaptation of a method introduced by \cite{BKT17} which provides this bijection. We can label each column of the bottom diagram as it follows: given $t\in [n]$, the $t$-th bottom column is
\begin{itemize}
\item $h$-related if there exists $i\in [k]$ such that $t = \alpha_{i}+i$.
\item $v$-related if there exists $j\in [n-k]$ such that $t = k+j-\alpha_{j}^{*}$. If $j$ exists then it must be unique.
\end{itemize}

There is a geometric interpretation in the HSYD: choose some bottom column and draw a 45-degree northwest line from the center of the first box in this column. If the line hits the last box in a top row, then the bottom column is $h$-related. Otherwise, it is $v$-related.

\begin{figure}[ht]
\centering
\includegraphics[scale=0.7]{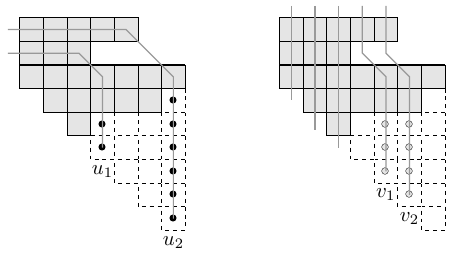}
\caption{On the left, we represent the $h$-related columns for $\alpha=(3,5)$ and $\lambda=(1,5,7)$. On the right, we represent the $v$-related columns. The vacant length of a column is the number of dot in the respective column}
\label{fig:related_columns}
\end{figure}

For instance, consider the permutation $w=2\,6|\overline{7}\,\overline{5}\,\overline{1}\,3\,4$ in Example \ref{explo:1}. Figure \ref{fig:related_columns} exhibits the corresponding double partition $\alpha=(3,5)$ and $\lambda=(1,5,7)$ for $w$ as well how the $h$-related columns are in bijection with rows in the top diagram (left) and how the $v$-related columns are in bijection with columns in the top diagram (right).\footnote{The terms $h$-related (short for ``horizontally'' related) and $v$-related (short for ``vertically'' related) were motivated by the bijection with the rows and column, resp.} This bijection is established in the next proposition.

\begin{prop}\label{prop:relatedcolumn}\ 
\begin{enumerate}
\item There is a bijection between $h$-related columns and rows in the top diagram.
\item There is a bijection between $v$-related columns and columns in the top diagram.
\item Any bottom column is either $h$-related or $v$-related;
\end{enumerate}
\end{prop}
\begin{proof}
For statement (1), if there are $i<i'$ such that $\alpha_{i}+i=\alpha_{i'}+i'$ then $0<i'-i = \alpha_{i}-\alpha_{i'}\leqslant 0$, which is impossible. Hence, two different rows in the top diagram are related to different $h$-related columns. Similarly, we have statement (2).

For statement (3), suppose that the $t$-th bottom column is, simultaneously, $h$-related and $v$-related. Then, there are $i$ and $j$ such that $t=\alpha_{i}+i=k+j-\alpha_{j}^{*}$, i.e., $\alpha_{i}-j=k-i-\alpha_{j}^{*}$. If $\alpha_{i}\geqslant j$ then $\alpha_{j}^{*}\leqslant k-i$ and $\alpha_{j}^{*}=\#\{l\tq \alpha_{l}\geqslant j\}\geqslant k-i+1$, a contradiction. On the other hand, if $\alpha_{i}< j$ then $\alpha_{j}^{*}>k-i$ and $\alpha_{j}^{*}=\#\{l\tq \alpha_{l}\geqslant j\}< k-i$, also a contradiction. Hence, no bottom column can be, simultaneously, $h$-related and $v$-related. Since the number of $h$-related and $v$-related bottom columns is respectively $k$ and $n-k$, a bottom column is either $h$-related or $v$-related.
\end{proof}

The \emph{vacant length of a bottom column} is the number of empty boxes below the boxes of $\lambda$ in the staircase $n\times n$ shape. Explicitly, the vacant length of the $j$-th bottom column is the number $j-\#\{i\tq \lambda_{i}+i> j\}$.

We may recover the permutation associated with such diagram by taking the vacant length of the $h$-related and $v$-related bottom columns. Namely, the permutation element for $\Lambda=(\alpha,\lambda)$ is defined by $w_{u,\lambda}$ in the Equation \eqref{eq:kgrass}, where $0<u_{1} < \cdots < u_{k}$ are the vacant length of the $h$-related columns, and $0<v_1< \cdots < v_{n-k-\ell(\lambda)}$ are the vacant length of the $v$-related columns. In Figure \ref{fig:related_columns}, the vacant length is the number of dots in the respective column. Then, $u_{1}= 2, u_{2}=6$ is the vacant length of the $h$-related columns and $v_{1}= 3, v_{2}=4$ is the vacant length of the $v$-related columns.

\subsection{HSYD's and covering types}\label{subsec:hsyd}

After the introduction of covering pairs and the definition of HSYD's, we now show how the diagrams illustrate covering relations.

We recall the following proposition about the covering relation in terms of double partitions.

\begin{prop}[\cite{LR19}, Proposition 4.2]\label{prop:doublepartitiontypes} Let $w,w'\in \weyl_{n}^{(k)}$. Denote by $\Lambda=(\alpha,\lambda)$ and $\Lambda'=(\alpha',\lambda')$ the associated double partitions of $w$ and $w'$, respectively. 
Then,
\begin{itemize}
\item $\typew(w,w')=$ B1 if, and only if, for every $i \in [k]$ and $j \in [r-1]$ we have
$\alpha'_{i}  = \alpha_{i}$ and $\lambda'_{j} = \lambda_{j+1}$.

\item $\typew(w,w')=$ B2 if, and only if, for every $i \in [k]$ and $j \in [r]$ we have
\begin{align*}
\alpha'_{i}  &= \alpha_{i}
& \mbox{ and } &&
\lambda'_{j}  &= 
\left\{
\begin{array}{cc}
\lambda_{j}-1 & \mbox{if } j=t \\ 
\lambda_{j} & \mbox{if } j\neq t
\end{array},
\right.
\end{align*}
for some $t\in [r]$.
\item $\typew(w,w')=$ B3 if, and only if, for every $i \in [k]$ and $j \in [r]$ we have
\begin{align*}
\alpha'_{i} &= 
\left\{
\begin{array}{cc}
\alpha_{i}-1 & \mbox{if } i=p \\ 
\alpha_{i} & \mbox{if } i\neq p
\end{array} 
\right.
& \mbox{ and } &&
\lambda'_{j} &= \lambda_{j},
\end{align*}
for some $p\in [k]$.

\item $\typew(w,w')=$ B4 if, and only if, for every $i \in [k]$ and $j \in [r]$ we have
\begin{align*}
\alpha'_{i} &= 
\left\{
\begin{array}{cc}
\alpha_{i}+x-1 & \mbox{if } i=p \\ 
\alpha_{i} & \mbox{if } i\neq p
\end{array} 
\right.
& \mbox{ and } &&
\lambda'_{j} &= 
\left\{
\begin{array}{cc}
\lambda_{j}-x & \mbox{if } j=t \\ 
\lambda_{j} &  \mbox{if } j\neq t
\end{array},
\right.
\end{align*}
for some $p\in [k]$ and $t\in[r]$.
\end{itemize}
\end{prop}

\begin{rem}By Proposition \ref{prop:doublepartitiontypes}, the covering relations in $\mathcal{P}(k,n)$ are almost the same as the covering relations for $\mathcal{D}_{k,n}$ and $\mathcal{SD}_{n}$, except in the case of a pair of type B4 when $x\neq 1$. This is exactly the step from the weak to strong Bruhat order. It also reflects the non-fully commutativeness of the class of $k$-Grassmannian permutations when $k\neq 0$.
\end{rem}

The results of the Proposition \ref{prop:doublepartitiontypes} may be understood in terms of the operation of removing boxes of the HSYD. There are two types of boxes that can be removed: the corners and the middle boxes. 

A \textit{corner} is a box of the diagram when removed produces a new diagram without any further operation. By Proposition \ref{prop:doublepartitiontypes}, we have that:
\begin{itemize}
\item if $\typew(w,w')=$ B1, the diagram of $\Lambda'$ is obtained by removing a corner in the diagonal of the bottom diagram;
\item if $\typew(w,w')=$ B2, the diagram of $\Lambda'$ is obtained from the diagram of $\Lambda$ by removing a corner of the bottom diagram that belongs to a $v$-related column;
\item if $\typew(w,w')=$ B3, the diagram of $\Lambda'$ is obtained from the diagram of $\Lambda$ by removing a corner of the top diagram.
\end{itemize}

A \textit{middle (bottom) box}, is a box (which is not a corner) of the bottom diagram when removed produces a new diagram after a displacement of boxes at its right. A middle bottom box is neither a corner, nor a diagonal; it lies in an $h$-related column; and all boxes to the right of it should belong to a $v$-related column. If $\typew(w,w')=$ B4, by Proposition \ref{prop:doublepartitiontypes}, the diagram of $\Lambda'$ is obtained from the diagram of $\Lambda$ by removing a box of the bottom diagram that belongs to an $h$-related column which is either a corner or a middle box, respectively, when $x=1$ or $x\neq 1$. Notice that if $x\neq 1$ then both partitions are changed which is a consequence of the movement of the $(x-1)$ boxes at the right of the removed box in the bottom diagram to the top diagram.

\begin{figure}[ht]
\centering
\includegraphics[scale=0.7]{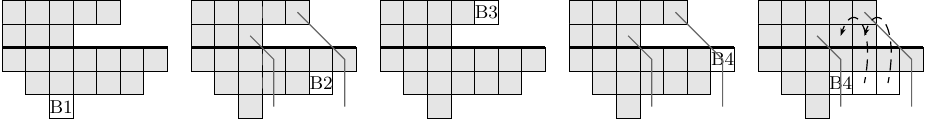}  
\caption{The five coverings of $w=2\,6|\overline{7}\, \overline{5}\, \overline{1}\, 3\, 4$ obtained by removing boxes in its HSYD according to the type of the pair.}
\label{fig:removing_boxes}
\end{figure}

In Example \ref{explo:1}, the permutation $w=2\,6|\overline{7}\,\overline{5}\,\overline{1}\,3\,4$ covers five different elements in $\weyl_{7}^{(k)}$. Figure \ref{fig:removing_boxes} illustrate these covering pairs in terms of removing corners or middle boxes. 

\subsection{Row-reading}

With the aid of the HSYD's, there is an easy way to get a reduced decomposition of a permutation $w_{\Lambda}$ using both the constructions for Grassmannians of type A and for maximal isotropic Grassmannians of type B given in \cite{IN08}. Let $w_{\alpha}$ be the row-reading given by the partition $\alpha$ and $w_{\lambda}$ be the row-reading given by the strict partition $\lambda$. More specifically,
\begin{itemize}
\item the row-reading of $\alpha$ is obtained by assigning a simple reflection consecutively to each box of the Young diagram from left to right and upwards, starting from $s_{1}$ in the bottom leftmost box. Then, $w_{\alpha}$ is the word obtained reading each row in the diagram from right to left, and the rows from bottom to top.

\item the row-reading of $\lambda$ is obtained by assigning a simple reflection consecutively to each box of the strict Young diagram from left to right, starting from $s_{0}$ in leftmost box of each row in the staircase diagram. Then, $w_{\lambda}$ is the word obtained reading each row in the diagram from right to left, and the rows from bottom to top.
\end{itemize}

Then, $w = w_{\lambda} \cdot w_{\alpha}$ is a reduced decomposition of such permutation called row-reading of $\Lambda$.

Hence, with respect to the reduced decompositions of $w$ and $w'$, the removing of a corner may be seen as getting $w'$ directly by deleting the corresponding reflection $s_{\ast}$ in the reduced decomposition of $w$ given by the row-reading. Now, the displacement of the boxes at the right  of a removed middle box occurs because the reduced expression of $w'$ is obtained by deleting some reflection $s_{\ast}$ in the reduced decomposition of $w$ followed by the application of some relations.

\begin{figure}[ht]
\centering
	\includegraphics[scale=0.7]{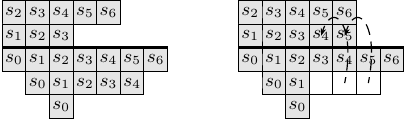}
	\caption{Row-reading of $w=2\,6|\overline{7}\,\overline{5}\,\overline{1}\,3\,4$ on the left and row-reading of $w'=5\,6|\overline{7}\,\overline{2}\,\overline{1}\,3\,4$ on the right}
	\label{fig:rowreading}
\end{figure}

For instance, the row-reading of $w=2\,6|\overline{7}\,\overline{5}\,\overline{1}\,3\,4$ is obtained by reading each row of Figure \ref{fig:rowreading} (left) backwards, from the lowest to the topmost row. In this case, we get $w = s_{0}\cdot s_{4}s_{3}s_{2}s_{1}s_{0}\cdot s_{6}s_{5}s_{4}s_{3}s_{2}s_{1}s_{0}\cdot s_{3}s_{2}s_{1}\cdot s_{6}s_{5}s_{4}s_{3}s_{2}$. Notice that $w$ covers $w'=5\,6|\overline{7}\,\overline{2}\,\overline{1}\,3\,4$, which corresponds to remove a middle box and resulting in $w' = s_{0}\cdot s_{4}s_{3} \widehat{s_{2}}s_{1}s_{0}\cdot s_{6}s_{5}s_{4}s_{3}s_{2}s_{1}s_{0}\cdot s_{3}s_{2}s_{1}\cdot s_{6}s_{5}s_{4}s_{3}s_{2}$. However, the row-reading of $w'$ comes from Figure \ref{fig:rowreading} (right) and is $w' = s_{0}\cdot s_{1}s_{0}\cdot s_{6}s_{5}s_{4}s_{3}s_{2}s_{1}s_{0}\cdot s_{5}s_{4}  s_{3}s_{2}s_{1}\cdot s_{6}s_{5}s_{4}s_{3}s_{2}$.

\begin{rem} This extends to any isotropic Grassmannians the reading map introduced in Section 4.2 of \cite{IN08} in the context of maximal isotropic Grassmannians.
\end{rem}

\section{Boundary map and integral homology}\label{sec:boundarymap}

In this section, after presenting the permutation model together with some combinatorial properties of the Weyl group of type B, we go into the details of the computation of the homology groups. 

We construct a chain complex $(\mathcal{C},\partial)$ whose homology is $H_{*}(\mathbb{F}_{(k)},\mathbb{Z})$, where $\mathbb{F}_{(k)}$ is the real flag manifold that is either the isotropic Grassmannian $\mathrm{IG}(n-k,2n)$ or the odd orthogonal Grassmannians $\mathrm{OG}(n-k,2n+1)$. We state below the general results found in \cite{RS19} adapted to our context.

Let $\mathcal{C}$ be the $\mathbb{Z}$-module freely generated by $\mathcal{S}_{w}$, for every element $w$ of the set of minimal representatives $\smash{\weyl_{n}^{(k)}}$. The boundary map $\partial \colon\mathcal{C} \rightarrow \mathcal{C}$ is defined by 
\begin{equation}\label{eq:boundcoeff}
\partial \mathcal{S}_{w}=\sum_{w^{\prime }}c(w,w^{\prime })\mathcal{S}_{w^{\prime }}
\end{equation}
for some coefficients $c(w,w^{\prime })\in \mathbb{Z}$. By \cite{RS19} Theorem 2.2, we know that $c(w,w^{\prime})$ is either $0$ or $\pm 2$. If $\ell(w)-\ell(w')\neq 1$ or $w'$ and $w$ are not comparable by the Bruhat-Chevalley order ``$\leqslant$'' then $c(w,w^{\prime })=0$. The non-zero coefficients may occur when $w'\leqslant w$ and $\ell(w)=\ell(w')+1$, i.e., for $w$ covering $w'$.

For $w \in \mathcal{W}^{(k)}_{n}$, define $\Pi_w=\Pi^+ \cap w\Pi^-$, the set of positive roots sent to negative roots by $w^{-1}$. Let $\sigma(w)$ be the sum of roots in $\Pi_{w}$, i.e.,
\begin{equation*}
\sigma (w)=\sum_{\beta \in \Pi _{w}} \beta.
\end{equation*}

\begin{prop}[\cite{RS19}, Proposition 2.7] Let $\gamma$ be the unique root (not necessarily simple) such that $w=s_{\gamma}w'$. Then
\begin{equation}\label{eq:kappa}
\sigma(w)-\sigma(w')=\kappa \cdot \gamma
\end{equation}
for some integer $\kappa=\kappa(w,w')$.
\end{prop}

\begin{thm}[\cite{RS19}, Thm. 2.8, \cite{Koc95}, Thm. 1.1.4]\label{thm:kocherlakota}
Suppose that $w$ covers $w'$. Then the coefficient $c(w,w')$ is given as follows:
\begin{equation*}
c(w,w')= \pm (1+(-1)^{\kappa})=
\left\{
\begin{array}{cl}
0 & \mbox{if } \kappa \mbox{ is odd,} \\ 
\pm2 & \mbox{if } \kappa \mbox{ is even.}
\end{array}
\right.
\end{equation*}
The signs on $c(w,w')$ can be chosen so that $\partial^{2}=0$ and the homology of $(\mathcal{C},\partial)$ is the integral homology of $\mathbb{F}_{(k)}$.
\end{thm}
\begin{rem} The formula for $c(w,w^{\prime})$ obtained in \cite{RS19} offers a choice of the signs defined in terms of the reduced decompositions for the Weyl group elements.
\end{rem}

Our objective is to use Theorem \ref{thm:kocherlakota} to compute the boundary coefficients of the isotropic and odd orthogonal Grassmannians, which consists in the computation of $\sigma(w)-\sigma(w')$ for $w$ covering $w'$, where $\sigma(w)=\sum_{\beta\in\Pi_{w}}\beta$ is the sum of all roots of $\Pi_{w}$. This process will provide us both $\kappa$ and $\gamma$. 

Our main strategy is based on a bijective correspondence between the roots of $\Pi_w$ with the half-shifted Young diagram of $w$ as we explain now.

\subsection{Inversions}\label{subsec:inversions}

The number of boxes of the HSYD of $w$ gives its length with $|\alpha|$  inversions of the first $k$ entries and $|\lambda|$ inversions corresponding to negative part. In this section, the HSYD's are used to describe $\Pi_w$ as the union of two distinguished sets given by the top and bottom diagrams of $w$ parametrizing the inversions of the corresponding partitions $\alpha$ and $\lambda$ respectively. For each $w\in \weyl_{n}^{(k)}$, we define
\begin{align}
\invp(w) & = \{(i,j)\in [n]^{2} \mid i<j \mbox{ and } w(i)>w(j)\} \label{eq:rootspos}\\
\invn(w) &= \{(i,j)\in [n]^{2} \mid i\leqslant j \mbox{ and } -w(i)>w(j)\}. \label{eq:rootsneg}
\end{align}

It follows that $|\invp(w)|+|\invn(w)| = |\alpha|+|\lambda|=\ell(w)$ by Equation \eqref{eq:length} . 

If $w=s_{i_{1}}\cdots s_{i_{\ell}}$ is a reduced decomposition of $w$ then it is known that $|\Pi_{w}|=\ell(w)$ and  $\Pi_{w} = \{a_{i_{1}}, \ s_{i_{1}}(a_{i_{2}}), \ s_{i_{1}}s_{i_{2}}(a_{i_{3}}),\dots, \ s_{i_{1}}\!\cdots s_{i_{\ell-1}}(a_{i_{\ell}})\}$. Consider $\beta^{+}_{i,j}$ and $\beta^{-}_{i,j}$ defined, resp., for inversions $(i,j)$ in $\invp(w)$ and $\invn(w)$ as follows
\begin{itemize}
\item For $(i,j)\in \invp(w)$,
\begin{equation*}
\beta^{+}_{i,j}(w) =\varepsilon_{w(i)} -\varepsilon_{w(j)}
\end{equation*}

\item For $(i,j)\in \invn(w)$,
\begin{equation*}
\beta^{-}_{i,j}(w) =
\left\{
\begin{array}{cl}
-\varepsilon_{{w(i)}}-\varepsilon_{w(j)}  & \mbox{, if type C},\\
2^{-\delta_{ij}}(-\varepsilon_{{w(i)}}-\varepsilon_{w(j)}) & \mbox{, if type B},
\end{array}
\right.
\end{equation*}
where $\delta_{ij}$ is the Kronecker delta.
\end{itemize}

\begin{prop}\label{prop:inversions}
Let $w\in \weyl_{n}^{(k)}$.
Then,
\begin{enumerate}
\item $\beta^{+}_{i,j}$ is a positive root in $\Pi_{w}$ for every $(i,j)\in\invp(w)$;
\item $\beta^{-}_{i,j}$ is a positive root in $\Pi_{w}$ for every $(i,j)\in\invn(w)$.
\end{enumerate}

Moreover, $\Pi_{w}$ is the disjoint union of $\beta^{+}$ and $\beta^{-}$, i.e., $\Pi_{w}= \{\beta^{+}_{i,j}\tq (i,j)\in\invp(w)\} \cup \{\beta^{-}_{i,j}\tq (i,j)\in\invn(w)\}$.
\end{prop}
\begin{proof}
For $(i,j)\in\invp(w)$ with $i<j$, the root $\varepsilon_{i}-\varepsilon_{j}$ is a negative root, whereas $w(\varepsilon_{i}-\varepsilon_{j})=\varepsilon_{w(i)}-\varepsilon_{w(j)}=\beta^{+}_{i,j}$ is a positive root. Then, $\beta^{+}_{i,j}\in \Pi_{w}$. For $(i,j)\in\invn(w)$ such that either $G$ a symplectic group (type C) with $i\leqslant j$ or $G$ an odd orthogonal group (type B) with $i<j$, the root $-\varepsilon_{i}-\varepsilon_{j}$ is a negative root, whereas $w(-\varepsilon_{i}-\varepsilon_{j})=-\varepsilon_{w(i)}-\varepsilon_{w(j)}=\beta^{-}_{i,j}$ is a positive root. Then, $\beta^{-}_{i,j}\in \Pi_{w}$. For $(i,i)\in\invn(w)$ such that $G$ an odd orthogonal group (type B), the root $-\varepsilon_{i}$ is a negative root, whereas $w(-\varepsilon_{i})=-\varepsilon_{w(i)}=\frac{1}{2}(-\varepsilon_{w(i)}-\varepsilon_{w(i)})=\beta^{-}_{i,j}$ is a positive root. Then, $\beta^{-}_{i,j}\in \Pi_{w}$.

Since $|\Pi_{w}|=|\invp(w)|+|\invn(w)|$ and, furthermore, all roots $\beta^{+}$ and $\beta^{-}$ are different to each other, the sets $\Pi_{w}$ and $\{\beta^{+}_{i,j}\tq (i,j)\in\invp(w)\} \cup \{\beta^{-}_{i,j}\tq (i,j)\in\invn(w)\}$ should be equal.
\end{proof}

We have a different description of the set $\invp(w)$ which will be useful.

\begin{prop}\label{prop:invprop_pos}
Given $w\in\weyl^{(k)}_{n}$, we have that \ $\invp(w) = \{(i,j) \tq i\in [k] \mbox{ and } j\in [k+1, k+\alpha_{i}]\}$.
\end{prop}
\begin{proof}
Recall that $w\in \weyl_{n}^{(k)}$ is given by $w = u_1\, \cdots\, u_k | \overline{\lambda_r}\, \cdots\, \overline{\lambda_1}\, v_{1}\, \cdots \, v_{n-k-r}$ and it satisfies Equations \eqref{eq:kgrass_conds}. 

Consider $I_{1}=\{(i,j) \tq  i\in [k] \mbox{ and } j\in[k+1, k+\alpha_{i}]\}$. Notice firstly that $|I_1|=\sum_{i=1}^{k} |\alpha_i|=|\alpha|=|\invp(w)|$ so that we only need to prove $I_{1}\subset \invp(w)$. Suppose that $(i,j)\in I_{1}$. Since $\alpha_{i}\geqslant r$, we can split in two cases: if $j\in [k+1,k+r]$ then $w(i)=u_{i}$ and $w(j)=\overline{\lambda}_{k+r-i+1}$, which clearly implies that $w(i)>w(j)$; if $j\in [k+r+1, k+\alpha_{i}]$ then, by Equation \eqref{eq:alpha2}, we have that $1\leqslant j-k-r \leqslant \#\{v_{l}\tq v_{l}<u_{i}\}$, i.e.,$w(j)=v_{j-k-r}<u_{i}=w(i)$. 
\end{proof}

We now seek a similar description for $\invn(w)$. However, the characterization does not follow so directly as that made for $\invp(w)$.. The strategy here is based on a relationship between $w$ and its corresponding Lagrangian permutation $\wmax$ throughout a ``translation'' map.   

Let $w\in\weyl_{n}^{(k)}$ be a permutation in the form \eqref{eq:kgrass}. The element $\wmax \in \weyl_{n}^{(0)}$ is obtained by reordering the entries in $w$ such that $\wmax = | \overline{\lambda_r}\, \cdots\, \overline{\lambda_1}\, \widetilde{v}_{1}\, \cdots \, \widetilde{v}_{n-r}$, where $0< \widetilde v_{1}< \cdots < \widetilde v_{n-r}$ is the reordering of the entries $u$'s and $v$'s of $w$. Next lemma allows us to keep track the position of the $w(i)$'s after this process. Indeed, if $w(i)=j$ then $w^{-1}(j)=i$, so that the position of the value $j$ in the permutation $w$ is equal to $w^{-1}(j)$. 

\begin{lem}\label{lema:lagran_perm}
Let $w\in\weyl_{n}^{(k)}$. The position of the entry $w(i)$,  for $ i \in [n]$, in the Lagrangian permutation $\wmax\in\weyl_{n}^{(0)}$ is given by the formula
\begin{equation*}
\wmax^{-1}(w(i)) =
\left\{
\begin{array}{cl}
i+\alpha_{i} & \mbox{if } 1\leqslant i \leqslant k;\\ 
i-\alpha_{i-k}^{*} & \mbox{if } k+1\leqslant i \leqslant n.
\end{array} 
\right.
\end{equation*}

\end{lem}
\begin{proof}
If $1\leqslant i\leqslant k$ then $\wmax^{-1}(w(i))= i+r+\#\{v_{j}\tq v_{j}<u_{i}\}=i+\alpha_{i}$, by Equation \eqref{eq:alpha2}.
If $k+1\leqslant i \leqslant k+r$ then, by Lemma \ref{lema:conjugationExplict}, we have $\wmax^{-1}(w(i))=i-k=i-\alpha_{i-k}^{*}$.
If $k+r+1\leqslant i \leqslant n$ then, by Lemma \ref{lema:conjugationExplict}, we have $\wmax^{-1}(w(i))=i-\#\{u_{j} \tq u_{j}>v_{i-k-r}\}=i-\mu_{i-k-r}^{*}=i-\alpha_{i-k}^{*}$.
\end{proof}

The point is that we may find a description for $\invn(\wmax)$, the set of negative inversions of the Lagrangian permutation $\wmax\in\weyl_{n}^{(0)}$ associated with $w$.

\begin{prop}\label{prop:invprop_neg}
Given $w\in\weyl^{(k)}_{n}$, we have that \ $\invn(\wmax) = \{(i,j) \tq  i\in [r] \mbox{ and } j\in [i, i-1+\lambda_{r-i+1}]\}$.
\end{prop}
\begin{proof}
Consider $I_{2} =\{(i,j) \tq i\in [r] \mbox{ and } j\in [i, i-1+\lambda_{r-i+1}]\}$. Observe that $\lambda_{l}\geqslant l$ for every $l \in [r]$ because $\lambda$ is a strict partition. Then, $i-1+\lambda_{r-i+1}\geqslant r$. Notice also that $|I_2|=\sum_{i=1}^{r} |\lambda_i|=\ell(\wmax)=|\invn(\wmax)|$ so that we only need to prove $I_{2}\subset \invn(\wmax)$. Suppose that $(i,j)\in I_{2}$ with $i\leqslant j$. If $j \in [r]$ then $w(i)=\overline{\lambda}_{r-i+1}$ and $w(j)=\overline{\lambda}_{r-j+1}$ so that $-w(i)>w(j)$. If $j\in [r+1, i-1+ \lambda_{r-i+1}]$ then $w(j)=\widetilde{v}_{j-r} \leqslant \widetilde{v}_{\lambda_{r-i+1}-(r-i+1)}$. Since $\widetilde{v}_{\lambda_{l}-l}<\lambda_{l}$ for every $l\in [r]$, we conclude that $w_{j}<\lambda_{r-i+1}=-w(i)$.
\end{proof}

Finally, a description of $\invn(\wmax)$ is quite enough for us since there is a bijection between $\invn(w)$ and $\invn(\wmax)$, as shown in the next lemma, which follows directly from the definition of inversions in $\invn(w)$ and by Lemma \ref{lema:lagran_perm}.

\begin{lem}\label{lema:inv_neg}
Given $w\in\weyl_{n}^{(k)}$, then $(i,j)\in\invn(w)$ if, and only if, $(i_{0},j_{0})\in\invn(\wmax)$, where $w(i)=\wmax(i_{0})$ and $w(j)=\wmax(j_{0})$. Moreover, $\beta_{i,j}^{-}(w)=\beta_{i_0,j_0}^{-}(\wmax)$.
\end{lem}

\begin{figure}[ht]
	\centering
	\includegraphics[scale=0.7]{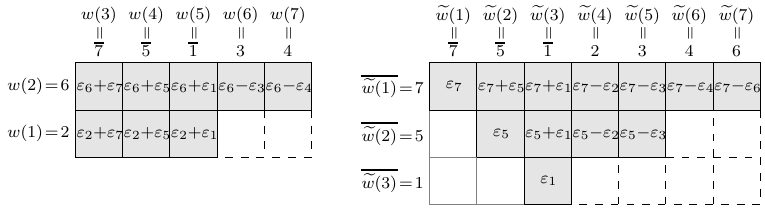}
	\caption{The roots of $\Pi_w$ inside the HSYD of $w=2\,6|\overline{7}\,\overline{5}\,\overline{1}\,3\,4 \in \mathrm{OG}(5,15)$}
	\label{fig:roots}
\end{figure}

In general, the desired correspondence between the $\Pi_w$ and the HSYD of $w$ according to results above is done as follows:
\begin{itemize}
\item Top diagram: for each $(i,j) \in \invp(w)$ with $i\in [k]$ and $j\in [k+1, k+\alpha_i]$, we place the root $\beta^{+}_{i,j} \in \Pi_w$ at the position $(i,j-k)$ in the corresponding top diagram $\mathcal{D}_{\alpha}$. 
\item Bottom diagram: for each $(i,j) \in \invn(\wmax)$ with $i\in [r]$ and $j\in [r+1, r+\lambda_i]$, we place the root $\beta^{-}_{i,j} \in \Pi_w$ at the position $(i,j)$ in the corresponding bottom diagram $\mathcal{SD}_{\lambda}$. 
\end{itemize}
Figure \ref{fig:roots} illustrate this procedure for $w=2\,6|\overline{7}\,\overline{5}\,\overline{1}\,3\,4 \in \mathrm{OG}(5,15)$.

\begin{rem} The Propositions \ref{prop:inversions}, \ref{prop:invprop_pos}, \ref{prop:invprop_neg} and Lemma \ref{lema:inv_neg} generalizes the Lemma 10 of \cite{IN08} in types B and C for any type of Grassmannian beyond the maximal ones. Notice that in the Figure \ref{fig:roots} we have both top and bottom diagrams.
\end{rem}

\subsection{Results}

The integral homology $H_{*}(\mathbb{F}_{(k)},\mathbb{Z})$ can be computed after we determine the coefficients of the boundary map according to the formula given in the Theorem \ref{thm:kocherlakota} and make a choice of signs. We provide an explicit expression of $\kappa$ in terms of the type of pairs.  In the sequence, we also present a visual interpretation to the statements by means of the HSYD's.

Recall that the integers $P$, $T$, and $Q$ correspond to the positions where $w$ changes when compared to $w'$ as in Equation \ref{eq:rightaction}.

\begin{thm}\label{thm:mainthm}
Let $w,w'$ be in $\weyl^{(k)}_{n}$ such that $w$ covers $w'$. Then, $\kappa$ as in Equation \eqref{eq:kappa} depend on the type of the pair $w,w'$ as follows:
\begin{itemize}
\item For $\mathrm{IG}(n-k,2n)$:
\begin{equation*}
\kappa(w,w') =
\left\{
\begin{array}{cc}
T & \mbox{if } \type(w,w') = B1 \\ 
T+Q & \mbox{if } \type(w,w') = B2 \\
Q-P & \mbox{if } \type(w,w') = B3 \\
P+T & \mbox{if } \type(w,w') = B4 \\
\end{array} 
\right.
\end{equation*}

\item For $\mathrm{OG}(n-k,2n+1)$:
\begin{equation*}
\kappa(w,w') =
\left\{
\begin{array}{cc}
2T-1 & \mbox{if } \type(w,w') = B1 \\ 
T+Q-1 & \mbox{if } \type(w,w') = B2 \\
Q-P & \mbox{if } \type(w,w') = B3 \\
P+T-1 & \mbox{if } \type(w,w') = B4 \\
\end{array} 
\right.
\end{equation*}

\end{itemize}
\end{thm}
\begin{rem} The determination of coefficients for classical (real) Grassmannians $\mathrm{Gr}(k,n)$ -- those of type A -- occurs as a particular case when $w,w'$ is of type B3. The Schubert varieties are parametrized by permutations of the symmetric group $S_n$ with a descent at $k$-th position. They correspond to the $\alpha$ partition and to poset $\mathcal{D}_{k,n}$ of the Young diagrams. Indeed, Theorem \ref{thm:mainthm} shows that if $w,w'$ is a covering pair and $w=w'\cdot (P,Q)$ then $\kappa(w,w')=Q-P$. 
%The interpretation in the Young diagram is given by the top diagram of the figure \ref{fig:B3_coeff}.
\end{rem}
The proof will be postponed to Section \ref{sec:proofs} (Propositions \ref{prop:mainB1}, \ref{prop:mainB2}, \ref{prop:mainB3} and \ref{prop:mainB4}).

We may also visualize the corresponding $\kappa$'s for each pair throughout the diagrams. %(Figure \ref{fig:coef_hom}).
We start by filling in each diagram of $w$ and $w'$ with the corresponding inversions of $\Pi_w$ and $\Pi_{w'}$, remembering that the diagram of $w'$ is obtained from the diagram of $w$ by removing either a corner or a middle box. The distinction between these diagrams may be displayed inside the diagram of $w$ as follows:
\begin{itemize}
\item fill in the removed box (r.b.) with $1$;
\item don't fill the boxes if the corresponding roots of $\Pi_w$ and $\Pi_{w'}$ are the same; 
\item fill in with either $\pm 1$ or $2$ (as it will be clarified below) if the corresponding roots of $\Pi_w$ and $\Pi_{w'}$ are different. %if they are different, then we put the corresponding contribution, either $\pm 1$ or $2$, inside the box. We explain below 
\end{itemize}

As a consequence of computations made in the Section \ref{sec:proofs}, it follows that $\kappa(w,w')$ is given as the sum of the numbers in the diagram. Since in each pair we have only two values (positions) changed in the permutation, the difference between the roots occurs along specific rows and columns depending on the type of the pair as we now explain.
 
\begin{figure}[ht]

\begin{subfigure}{.45\textwidth}
	\centering
	\includegraphics[scale=0.75]{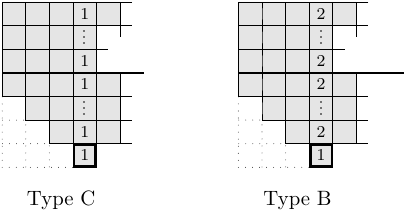}
	\caption{$\typew(w,w')=$ B1.}
	\label{fig:B1_coeff}
\end{subfigure}
\begin{subfigure}{.45\textwidth}
	\centering
	\includegraphics[scale=0.75]{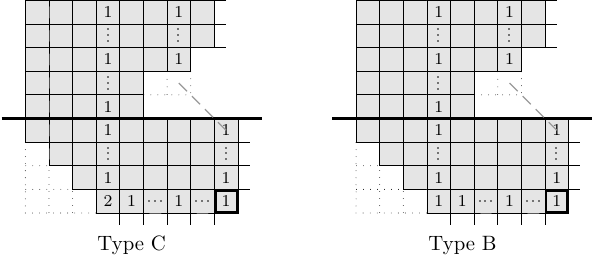}
	\caption{$\typew(w,w')=$ B2.}
	\label{fig:B2_coeff}
\end{subfigure}
\vspace*{1em}

\begin{subfigure}{.45\textwidth}
	\centering
	\includegraphics[scale=0.75]{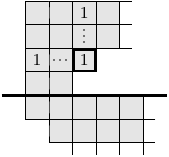}
	\caption{$\typew(w,w')=$ B3.}
	\label{fig:B3_coeff}
\end{subfigure}
\begin{subfigure}{.45\textwidth}
	\centering
	\includegraphics[scale=0.75]{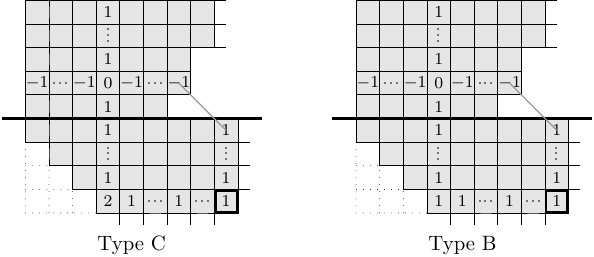}
	\caption{$\typew(w,w')=$ B4.}
	\label{fig:B4_coeff}
\end{subfigure}
\caption{Filling in the diagrams for each type of pair.}
\end{figure}

\begin{itemize}
\item If $\type(w,w')=$ B1, then we fill in the column above the r.b. in the diagonal of the bottom diagram according to Figure \ref{fig:B1_coeff}.
\item If $\type(w,w')=$ B2, remember the r.b. belongs to a $v$-related column in the bottom diagram. We fill in the hook defined by the diagonal box contained in the row of the r.b, and the column above the r.b. in bottom diagram together with its corresponding column in the top diagram according to Figure \ref{fig:B2_coeff}.
\item If $\type(w,w')=$ B3, remember the r.b. belongs to the top diagram. We fill in the boxes to the left and above the r.b. according to Figure \ref{fig:B3_coeff}.
\item If $\type(w,w')=$ B4, remember the r.b. belongs to anfir $h$-related column in the bottom diagram. We fill in the hook defined by the diagonal box contained in the row of the r.b, and the column above the r.b. in bottom diagram together with its corresponding line in the top diagram according to Figure \ref{fig:B4_coeff}.
\end{itemize}

\begin{ex}\label{explo:2} Let $w=2\,6|\overline{7}\,\overline{5}\,\overline{1}\,3\,4 \in \mathrm{OG}(5,15)$, where $n=7$, $k=2$. There are five covering pairs $w,w_i'$, $i=1,\ldots,5$, according to Example \ref{explo:1}.
We can compute $\kappa$ of each pair $w,w'_{i}$ as it follows:
\begin{align*}
\mbox{ (B1) } & w, w'_{1}: T=5 \mbox{ and } \kappa = 2T-1=9;\\
\mbox{ (B2) } & w, w'_{2}: (T,Q)=(4,7)  \mbox{ and } \kappa = T+Q-1=10 \\
\mbox{ (B3) } & w, w'_{3}: (P,Q)=(2,7) \mbox{ and } \kappa = Q-P=5; \\
\mbox{ (B4) } & w, w'_{4}: (P,T)=(2,3) \mbox{ and } \kappa = P+T-1= 4;\\
\mbox{ (B4) } & w, w'_{5}: (P,T)=(1,4) \mbox{ and } \kappa = P+T-1= 4.
\end{align*}

Therefore, by Theorem \ref{thm:mainthm}, $\partial{\schub}_{2\,6|\overline{7}\,\overline{5}\,\overline{1}\,3\,4} = \pm 2 \schub_{2\,4|\overline{7}\,\overline{5}\,\overline{1}\,3\,6} \pm 2 \schub_{2\,7|\overline{6}\,\overline{5}\,\overline{1}\,3\,4} \pm 2 \schub_{5\,6|\overline{7}\,\overline{2}\,\overline{1}\,3\,4}$.
\begin{figure}[ht]
	\centering
	\includegraphics[scale=0.7]{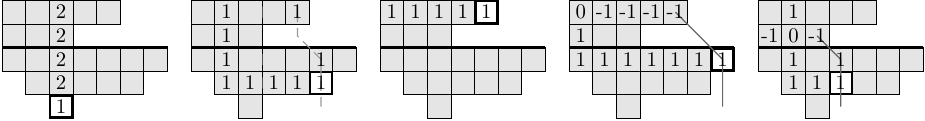}
	\caption{The sum of the number in each diagram gives the value of $\kappa(w,w_{i}')$.}
	\label{fig:coef_hom}
\end{figure}
\end{ex}

The comparison between the formulas for $\kappa(w,w')$ given in Theorem \ref{thm:mainthm} and the description of the covering relations in terms of the right action by Equation \eqref{eq:rightaction} had revealed a hidden phenomenon which will be generalized in a forthcoming paper.

Let us define the coroot of a root $\delta\in\Pi$ by $\delta^{\vee}=\dfrac{2\delta}{\langle \delta,\delta\rangle}$. The set of coroots $\Pi^{\ast}$ is also a root system which is called the dual root system. It is important to point out that root systems of type B and C are dual to each other.

If $\delta \in \Pi$ is given in terms of the system $\Sigma$ of simple roots as $\delta = \sum_{\xi \in \Sigma} d_{\xi}\xi  $, the height of the root $\delta$ is defined by the sum $\mathrm{ht}(\delta)=\sum_{\xi \in \Sigma} d_{\xi}$.

\begin{thm}\label{thm:dualheight}
Let $w,w'$ be in $\weyl^{(k)}_{n}$ such that $w$ covers $w'$ and $\delta$ be the root for which $w = w'\cdot s_{\delta}$. Then 
\begin{equation*}
\kappa(w,w')= \mathrm{ht}(\delta^{\vee}).
\end{equation*}
\end{thm}
\begin{proof}
Recall that the simple roots of type B are defined by $a_{0}=\varepsilon_{1}$ and $a_{i}=\varepsilon_{i+1}-\varepsilon_{i}$ for $1\leqslant i < n$. The positive roots $\Pi^+$ of type B are given by $\varepsilon_j$ ($j\geq 1$), $\varepsilon_i-\varepsilon_j$ ($i>j$) and $\varepsilon_i+\varepsilon_j$ ($i>j)$. The height of the corresponding roots and coroots 
\begin{align*}
\mathrm{ht}(\varepsilon_j) &=j , \mbox{ for } j\geqslant 1\\
\mathrm{ht}(\varepsilon_i-\varepsilon_j) &= i-j , \mbox{ for } i>j\\
\mathrm{ht}(\varepsilon_i+\varepsilon_j) &= i+j , \mbox{ for } i>j
\end{align*}

Now, recall that the simple roots of type C are defined by $a_{0}=2\varepsilon_{1}$ and $a_{i}=\varepsilon_{i+1}-\varepsilon_{i}$ for $1\leqslant i < n$. The positive roots $\Pi^+$ of type B are given by $2\varepsilon_j$ ($j\geq 1$), $\varepsilon_i-\varepsilon_j$ ($i>j$) and $\varepsilon_i+\varepsilon_j$ ($i>j)$. The height of the corresponding roots and coroots 
\begin{align*}
\mathrm{ht}(2\varepsilon_j) &=2j-1 , \mbox{ for } j\geqslant 1\\
\mathrm{ht}(\varepsilon_i-\varepsilon_j) &= i-j , \mbox{ for } i>j\\
\mathrm{ht}(\varepsilon_i+\varepsilon_j) &= i+j-1 , \mbox{ for } i>j
\end{align*}

Both root systems are dual to each other by the relation:
\begin{itemize}
\item For the type B root $\varepsilon_j$, its coroot is the type C root $2\varepsilon_j$;
\item For the type B root $\varepsilon_i-\varepsilon_j$, its coroot is the type C root $\varepsilon_i-\varepsilon_j$;
\item For the type B root $\varepsilon_i+\varepsilon_j$, its coroot is the type C root $\varepsilon_i+\varepsilon_j$.
\end{itemize}

Considering $w=w'\cdot s_{\delta}$, Equation \eqref{eq:rightaction} and Theorem \ref{thm:mainthm}, we have the following:
\begin{itemize}
\item For type C:
\begin{itemize}
\item If $\typew(w,w')=$ B1 then $\delta = 2\varepsilon_T$ and $\kappa(w,w')= T = \mathrm{ht}(\delta^{\vee})$;
\item If $\typew(w,w')=$ B2 then $\delta = \varepsilon_Q + \varepsilon_T$ and $\kappa(w,w')= T+Q = \mathrm{ht}(\delta^{\vee})$;
\item If $\typew(w,w')=$ B3 then $\delta = \varepsilon_Q - \varepsilon_P$ and $\kappa(w,w')= Q-P = \mathrm{ht}(\delta^{\vee})$;
\item If $\typew(w,w')=$ B4 then $\delta = \varepsilon_T + \varepsilon_P$ and $\kappa(w,w')= P+T = \mathrm{ht}(\delta^{\vee})$.
\end{itemize}
\item For type B:
\begin{itemize}
\item If $\typew(w,w')=$ B1 then $\delta = \varepsilon_T$ and $\kappa(w,w')= 2T-1 = \mathrm{ht}(\delta^{\vee})$;
\item If $\typew(w,w')=$ B2 then $\delta = \varepsilon_Q + \varepsilon_T$ and $\kappa(w,w')= T+Q-1 = \mathrm{ht}(\delta^{\vee})$;
\item If $\typew(w,w')=$ B3 then $\delta = \varepsilon_Q - \varepsilon_P$ and $\kappa(w,w')= Q-P = \mathrm{ht}(\delta^{\vee})$;
\item If $\typew(w,w')=$ B4 then $\delta = \varepsilon_T + \varepsilon_P$ and $\kappa(w,w')= P+T-1 = \mathrm{ht}(\delta^{\vee})$. \qedhere
\end{itemize}

\end{itemize}
\end{proof}

\subsection{Orientability}
We may provide general criteria of orientability.
\begin{prop}\label{prop:orient} \ 
\begin{enumerate}
\item $\mathrm{IG}(n-k,2n)$ is orientable if, and only if, $n-k$ is odd.

\item  $\mathrm{OG}(n,2n+1)$ is orientable for every $n$

\item $\mathrm{OG}(n-k,2n+1)$ is orientable if, and only if, $k>0$ and $n-k$ is even.
\end{enumerate}
\end{prop}
\begin{proof} Being orientable is equivalent to the boundary map in the top cell equals zero. 
The top cell $\schub_{w_{0}}$ corresponding to the longest element $w_{0} \in \weyl^{(k)}$ is given by Equation \eqref{eq:longestw}. If $k=0$ then the only possible choice is $w_{0}= | \overline{n\vphantom{1}}\,\overline{n-1}\, \cdots\, \overline{1}$ and $w'= | \overline{n\vphantom{1}}\,\overline{n-1}\, \cdots\, {1}$, which is a pair of type B1. In this case, $T=n$ and $P=Q=0$, implying that
\begin{equation*}
c(w_{0},w') =
\left\{
\begin{array}{cl}
\pm (1+(-1)^{n}) & \mbox{for } \mathrm{IG}(n,2n); \\ 
0 & \mbox{for } \mathrm{OG}(n,2n+1).
\end{array} 
\right.
\end{equation*}
Therefore, $\mathrm{IG}(n,2n)$ is orientable if, and only if, $n$ is odd, and $\mathrm{OG}(n,2n+1)$ is orientable for every $n$.
If $k>0$ then there is only one possible choice of $w'$ such that $w$ cover $w'$, namely, $w_{0} =1\,2\,\cdots\, (k-1)\, k| \overline{n\vphantom{1}}\,\overline{n-1}\, \cdots\, \overline{k+1}$ and $w' =1\,2\,\cdots\, (k-1)\, (k+1)| \overline{n\vphantom{1}}\,\overline{n-1}\, \cdots\, \overline{k}$, which is a pair of type B4. Then, $P=k$, $T=n$, and $Q=0$, implying that
\begin{equation*}
c(w_{0},w') =
\left\{
\begin{array}{cl}
\pm (1+(-1)^{k+n}) & \mbox{for } \mathrm{IG}(n-k,2n); \\ 
\pm (1-(-1)^{k+n}) & \mbox{for } \mathrm{OG}(n-k,2n+1).
\end{array} 
\right.
\end{equation*}
Therefore, $\mathrm{IG}(n-k,2n)$ is orientable if, and only if, $k+n\equiv n-k \mod 2$ is odd, and $\mathrm{OG}(n-k,2n+1)$ is orientable if, and only if, $k+n\equiv n-k \mod 2$ is even.
\end{proof}

\subsection{Duality}
Given $w\in \weyl^{(k)}_{n}$, define $w^{\vee}=w w_{0}$ the \emph{dual permutation of $w$}. If $w$ is written as in Equation \eqref{eq:kgrass}, the one-line notation of the dual permutation of $w$ is $w^{\vee} = u_1\,\cdots\, u_k | \overline{v_{n-k-r}}\,\cdots\, \overline{v_{1}}\, \lambda_{1}\,\cdots\, \lambda_{r}$.
The length of $w^{\vee}$ is $\ell(w^{\vee})=\ell(w_{0})-\ell(w)$. Next proposition states that the duality of a permutation also imply a duality over the covering pairs.

\begin{prop}[\cite{LR19}]\label{prop:dualtype} Let $w, w $ be permutations in $ \weyl^{(k)}_{n}$. Then, $w$ covers $w'$ if, and only if, $(w')^{\vee}$ covers $w^{\vee}$. Moreover,
\begin{enumerate}
\item $\typew(w,w')=$ B1 if, and only if, $\typew((w')^{\vee},w^{\vee})=$ B1;
\item $\typew(w,w')=$ B2 if, and only if, $\typew((w')^{\vee},w^{\vee})=$ B2;
\item $\typew(w,w')=$ B3 if, and only if, $\typew((w')^{\vee},w^{\vee})=$ B4;
\end{enumerate}
\end{prop}
The next proposition shows that $c((w')^{\vee},w^{\vee})$ can be obtained, according to $c(w,w')$ and the type of the pair $w,w'$.

\begin{prop}\label{prop:dualcoeff} Let $w,w'$ be in $\weyl^{(k)}_{n}$ such that $w$ covers $w'$. Then, $|c(w,w')|= |c((w')^{\vee},w^{\vee})|$ if, and only if, one of the following happens:
\begin{itemize}
\item For $\mathrm{IG}(n-k,2n)$: 
\begin{enumerate}
\item $\typew(w,w')=$ B1 and $n-k$ is odd;
\item $\typew(w,w')=$ B2;
\item $\typew(w,w')=$ B3 or B4, and $n-k$ is odd.
\end{enumerate}

\item For $\mathrm{OG}(n-k,2n+1)$:
\begin{enumerate}
\item $\typew(w,w')=$ B1 or B2;
\item $\typew(w,w')=$ B3 or B4, and $n-k$ is even.
\end{enumerate}

\end{itemize}
\end{prop}
\begin{proof}
Denote by $P^{\vee}$, $T^{\vee}$, and $Q^{\vee}$ the integers of Table \ref{tbl:ptq} for the dual pair $(w')^{\vee},w^{\vee}$. We need to describe such integers in terms of $P$, $T$, and $Q$, which should be considered case-by-case. Assume the flag manifold is $\mathrm{IG}(n-k,2n)$.

If $\typew(w,w')=$ B1 then $T^{\vee} =n-r+1=(n-k+1-2r)+T$. Thus, $(-1)^{T^{\vee}}=(-1)^{n-k+1} (-1)^{T}$ and we conclude that both coefficients are equal iff $n-k+1$ is even.

If $\typew(w,w')=$ B2 then $T^{\vee} =n-r-q+1=(n-k-2q-2r+1)+Q$ and $Q^{\vee} =n-r+t = (n-k -2r+2t-1) +T$. Thus, $(-1)^{T^{\vee}+Q^{\vee}}=(-1)^{T+Q}$ and we conclude that both coefficients are always equals.

If $\typew(w,w')=$ B3 then $P^{\vee} = p =P$ and $T^{\vee} =n-r-q+1=(n-k-2q-2r+1)+Q$. Thus, $(-1)^{P^{\vee}+T^{\vee}}=(-1)^{n-k+1} (-1)^{Q-P}$ and we conclude that both coefficients are equal iff $n-k+1$ is even.

If $\typew(w,w')=$ B4 then $P^{\vee} = p = P$ and $Q^{\vee} =n-r+t = (n-k -2r+2t-1) +T$. Thus, $(-1)^{Q^{\vee}-P^{\vee}}=(-1)^{n-k+1} (-1)^{P+T}$ and we conclude that both coefficients are equal iff $n-k+1$ is even.

The proof for $\mathrm{OG}(n-k,2n+1)$ is analogous.
\end{proof}

We define the \emph{incidence graph of a Grassmannian} as the graph whose vertices are the permutations of $\weyl_{n}^{(k)}$ and the edges ``$\rightarrow$'' and ``$\Rightarrow$'' are given by the covering relations as follows: if $w$ covers $w'$ and $c(w,w')=0$ then $w \rightarrow w'$; if $w$ covers $w'$ and $c(w,w')=\pm 2$ then $w \Rightarrow w'$. 

%For instance, 

\begin{ex}\label{ex:incidencegraph} The incidence graph for the odd orthogonal Grassmannian $\mathrm{OG}(2,9)$ and isotropic Grassmannian $\mathrm{IG}(2,8)$ where $n=4$, $k=2$, is given in Figure \ref{fig:n4k2_incidence}. The duality of permutations and pairs can be seen as the symmetry through the horizontal dashed line.

\begin{figure}[ht]
	\centering
	\includegraphics[scale=.8]{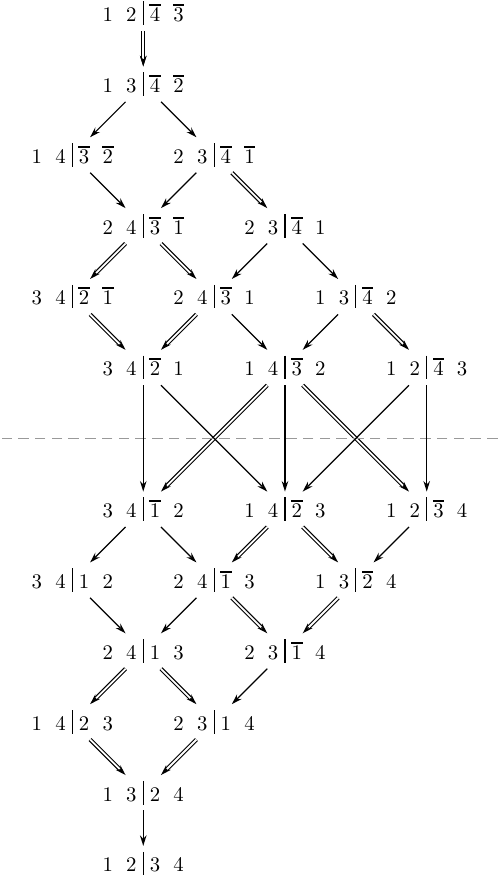}
	\hspace{2em}
	\includegraphics[scale=.8]{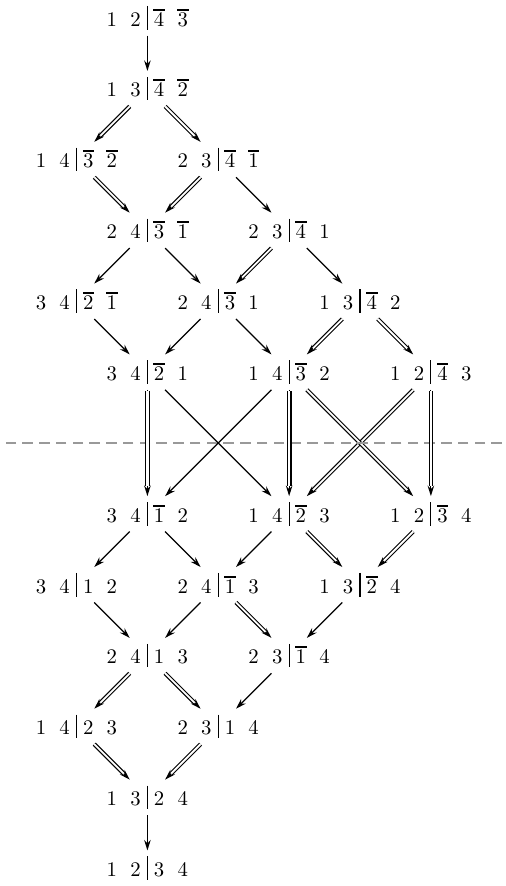}
	\caption{Incidence graph of $\mathrm{IG}(2,8)$ on the left and $\mathrm{OG}(2,9)$ on the right.}
	\label{fig:n4k2_incidence}
\end{figure}

%Therefore, their homology groups are
\begin{enumerate}
\item For $\mathbb{F}=\mathrm{IG}(2,8)$, by Proposition \ref{prop:dualcoeff}, only dual pairs $w,w'$ of type B2 satisfy $|c(w,w')| = |c((w')^{\vee},w^{\vee})|$. In the incidence graph in Figure \ref{fig:n4k2_incidence}, this means that the edge of a pair of type B1, B3, or B4 should be different from the edge associated with its dual. Therefore, the homology groups are
\begin{align*}
H_{11}(\mathbb{F},\Z) &= 0, &
H_{8}(\mathbb{F},\Z) &= \Z_{2}, &
H_{5}(\mathbb{F},\Z) &= \Z \oplus \Z_{2}, & 
H_{2}(\mathbb{F},\Z) &= \Z_{2}, \\
H_{10}(\mathbb{F},\Z) &= \Z_{2},&
H_{7}(\mathbb{F},\Z) &= \Z_{2}, &
H_{4}(\mathbb{F},\Z) &= \Z \oplus \Z_{2}, &
H_{1}(\mathbb{F},\Z) &= \Z_{2}, \\
H_{9}(\mathbb{F},\Z) &= \Z, &
H_{6}(\mathbb{F},\Z) &= (\Z_{2})^{2}, &
H_{3}(\mathbb{F},\Z) &= \Z_{2}, &
H_{0}(\mathbb{F},\Z) &= \Z.
\end{align*}

\item For $\mathbb{F}=\mathrm{OG}(2,9)$, by Proposition \ref{prop:dualcoeff}, the coefficients and the corresponding edges in Figure \ref{fig:n4k2_incidence} are the same for all dual pairs. Therefore, the homology groups are 
\begin{align*}
H_{11}(\mathbb{F},\Z) &= \Z, &
H_{8}(\mathbb{F},\Z) &= \Z_{2}, &
H_{5}(\mathbb{F},\Z) &= (\Z_{2})^{2}, & 
H_{2}(\mathbb{F},\Z) &= \Z_{2}, \\
H_{10}(\mathbb{F},\Z) &= 0,&
H_{7}(\mathbb{F},\Z) &= \Z \oplus \Z_{2}, &
H_{4}(\mathbb{F},\Z) &= \Z \oplus \Z_{2}, &
H_{1}(\mathbb{F},\Z) &= \Z_{2}, \\
H_{9}(\mathbb{F},\Z) &= \Z_{2}, &
H_{6}(\mathbb{F},\Z) &= \Z_{2}, &
H_{3}(\mathbb{F},\Z) &= \Z_{2}, &
H_{0}(\mathbb{F},\Z) &= \Z.
\end{align*}
\end{enumerate}

%The duality of permutations and pairs can be seen in Figure \ref{fig:n4k2_incidence} as the symmetry through the horizontal dashed line. Since $n-k=2$ is even, Proposition \ref{prop:dualcoeff} shows us that dual pairs $w,w'$ of type B1, B3, and B4 satisfy that $|c(w,w')|\neq |c((w')^{\vee},w^{\vee})|$. In the incidence graph in Figure \ref{fig:n4k2_incidence}, this means that the edge of a pair of type B1, B3, or B4 should be different from the edge associated with its dual.
\end{ex}

\subsection{Low dimensional homology groups}\label{subsec:lowdim}

Finally, as a by-product of our methods to compute the coefficients for the boundary map, it is not difficult to obtain results about the low-dimensional topology of these Grassmannians. %This result matches the one in del Barco-San Martin \cite{BS19}.

\begin{prop} \ 
\begin{enumerate}
\item If $\mathbb{F}=\mathrm{IG}(n-k,2n)$ then
\begin{align*}
H_{1}(\mathbb{F}, \Z) &=
\left\{
\begin{array}{cl}
\Z & \mbox{if } k=0, \\
\Z_{2} & \mbox{otherwise.}
\end{array} 
\right. ; &
H_{2}(\mathbb{F}, \Z) &=
\left\{
\begin{array}{cl}
0 & \mbox{if } n=2, k=1, \\
\Z_{2} & \mbox{otherwise.}
\end{array} 
\right.
\end{align*}
\item If $\mathbb{F}=\mathrm{OG}(n-k,2n+1)$ then
\begin{align*}
H_{1}(\mathbb{F}, \Z) &=
\left\{
\begin{array}{cl}
\Z & \mbox{if } n=1, k=0, \\
\Z & \mbox{if } n=2, k=1, \\
\Z_{2} & \mbox{otherwise.}
\end{array} 
\right. ; &
H_{2}(\mathbb{F}, \Z) &=
\left\{
\begin{array}{cl}
0 & \mbox{if } k=0, \\
\Z_{2} & \mbox{otherwise.}
\end{array} 
\right.
\end{align*}
\end{enumerate}
\end{prop}
\begin{proof}
To compute 1- and 2-homology, we only require to know the boundary maps $\partial_{1}$, $\partial_{2}$, and $\partial_{3}$, which depends on $k$ and $n$. Table \ref{tbl:coeff12homol} shows the incidence graph up to $3$-dimensional cells for different arrangements of $n$ and $k$. In the last graph, we denote $k_{a}=k+a$.

\begin{table}[ht]
\centering
\caption{Bruhat graphs to compute 1-, 2-homology}
\label{tbl:coeff12homol}
\begin{tabular}{c|c|c|c|c|c|c}
\hline 
& $n=1$& \multicolumn{2}{c|}{$n=2$} & \multicolumn{3}{c}{$n\geqslant 3$}\\ 
\hline 
& $k=0$ & $k=0$ & $k=1$ & \multicolumn{2}{c|}{$k=0$} & $k=1$ \\ 
\hline 
\begin{sideways}
Isotropic
\end{sideways} & 
\includegraphics[scale=.6]{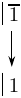} & 
\includegraphics[scale=.6]{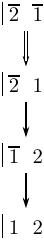} & 
\includegraphics[scale=.6]{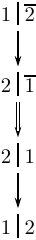} &
\multicolumn{2}{c|}{\includegraphics[scale=.6]{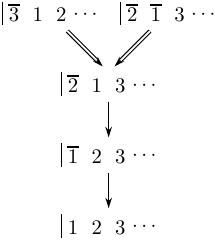}} &
\includegraphics[scale=.6]{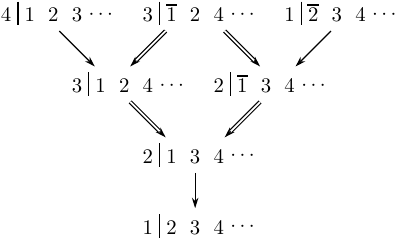} \\ 
\hline 
\begin{sideways}
Odd orthogonal
\end{sideways} &
\includegraphics[scale=.6]{figTable_Cn1k0.pdf} & 
\includegraphics[scale=.6]{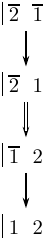} & 
\includegraphics[scale=.6]{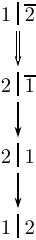} & 
\multicolumn{2}{c|}{\includegraphics[scale=.6]{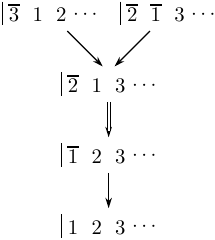}} &
\includegraphics[scale=.6]{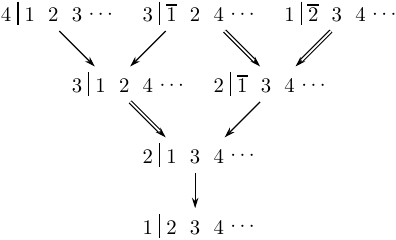}  \\ 
\hline \hline
&  \multicolumn{6}{c}{$n\geqslant 3$}\\ 
\hline 
& \multicolumn{4}{c}{$k=2$} & \multicolumn{2}{|c}{$k\geqslant 3$} \\ 
\hline 
\begin{sideways}
Any Grassm.
\end{sideways} & 
\multicolumn{4}{c}{\includegraphics[scale=.5]{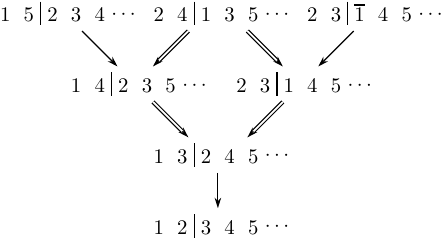}} &
\multicolumn{2}{|c}{\includegraphics[scale=.5]{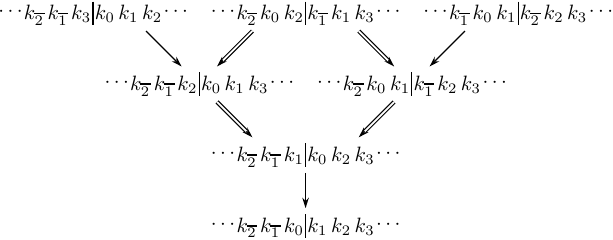}} \\ 
\hline 
\end{tabular} 
\end{table}

The 1-,2-homology groups can be computed using such diagrams.
\end{proof}

%\section{Proof of Theorem \ref{thm:mainthm}}\label{sec:proofs}
\section{Proof of the main results}\label{sec:proofs}

In this section, we compute $\sigma(w)-\sigma(w')$ for each type of covering pair $w,w'$. Remember the definition of the set of inversions of $w,w'$ as given in \eqref{eq:rootspos} and \eqref{eq:rootsneg}.

Let us denote $(\beta')^{+}_{i,j}=\beta^{+}_{i,j}(w')$ where $(i,j)\in\invp(w')$, $\widetilde\beta^{-}_{i,j}=\beta^{-}_{i,j}(\wmax)$ where $(i,j)\in\invn(\wmax)$, and $(\widetilde\beta')^{-}_{i,j}=\beta^{-}_{i,j}(\wprimemax)$ where $(i,j)\in\invn(\wprimemax)$. It will be useful to write $\sigma(w)-\sigma(w')=S^{+}+S^{-}$ with
\begin{equation*}
S^{+} =\quad \sum_{\mathclap{(i,j)\, \in \, \invp(w)}}\ \beta^{+}_{i,j} \quad -\quad\sum_{\mathclap{(i,j)\, \in \, \invp(w')}}\ (\beta')_{i,j}^{+}
\quad \mbox{ and } \quad
S^{-} =\quad \sum_{\mathclap{(i,j)\, \in \, \invn(\wmax)}}\ \widetilde\beta^{-}_{i,j} \quad -\quad\sum_{\mathclap{(i,j) \, \in \, \invn(\wprimemax)}}\ (\widetilde\beta')_{i,j}^{-}.
\end{equation*}

We can see in Theorem \ref{thm:mainthm} the $\kappa$ depends on the type of $G$. We can merge both types in a single formula after adopting the following notation: let $\isB$ be the variable that indicates whether the Grassmannian is either odd orthogonal or not, i.e.,
%$\isB= 0$ if $\mathbb{F}_{(k)}=\mathrm{IG}(n-k,2n)$, and $\isB= 1$ if $\mathbb{F}_{(k)}=\mathrm{OG}(n-k,2n+1)$
\begin{equation*}
\isB = 
\left\{
\begin{array}{cl}
1 & \mbox{ for } \mathrm{OG}(n-k,2n+1); \\ 
0 & \mbox{ for } \mathrm{IG}(n-k,2n).
\end{array} 
\right.
\end{equation*}

Then, $\kappa$ in Theorem \ref{thm:mainthm} can be given as following:
\begin{itemize}
\item Type B1: $\kappa = (1+\isB) T -\isB$;
\item Type B2: $\kappa = T+Q -\isB$;
\item Type B3: $\kappa = Q-P$;
\item Type B4: $\kappa = P+T-\isB$;
\end{itemize}

\begin{prop}\label{prop:mainB1}
Let $w,w'$ be in $\weyl^{(k)}_{n}$ such that $\typew(w,w')=$ B1. Denote $T=w^{-1}(\overline{\lambda_{1}})=k+r$. Then, $\kappa = (\isB+1)T -\isB$ and $\gamma = 2^{1-\isB}\ \varepsilon_{1}$.
\end{prop}
\begin{proof} Observe that $w(T) = \overline{\lambda_{1}} =\overline{1}$, $w'(T) = \lambda_{1} =1$, and $w(i)=w'(i)$ for every $i$ different from $T$. Recall that $\alpha_{i}'=\alpha_{i}$ for every $i$ (Proposition \eqref{prop:doublepartitiontypes}). Then, by Proposition \ref{prop:invprop_pos}, $\invp(w)= \invp(w')= \{(i,j) \tq i\in [k] \mbox{ and } j\in [k+1, k+\alpha_{i}]\}$.

Clearly, $T\in [k+1,k+\alpha_{i}]$ for every $i\in [k]$. For $(i,j)\in\invp(w)$, 
\begin{equation*}
\beta^{+}_{i,j}  = %\varepsilon_{w(i)}-\varepsilon_{w(j)} =
\left\{
\begin{array}{cl}
\varepsilon_{w(i)}+\varepsilon_{1} & \mbox{if } j=T\\ 
 \varepsilon_{w(i)}-\varepsilon_{w(j)} & \mbox{if } j\neq T
\end{array} 
\right. ; \quad
(\beta')^{+}_{i,j}  = %\varepsilon_{w'(i)}-\varepsilon_{w'(j)} =
\left\{
\begin{array}{cl}
\varepsilon_{w(i)}-\varepsilon_{1} & \mbox {if } j=T\\ 
 \varepsilon_{w(i)}-\varepsilon_{w(j)} & \mbox{if } j\neq T
\end{array} 
\right.
\end{equation*}

Then, $S^{+}$ can be rearranged as
\begin{equation*}
S^{+}=\sum_{i=1}^{k}\left( \sum_{j=k+1}^{k+\alpha_{i}} ( \beta^{+}_{i,j}-(\beta')^{+}_{i,j} )\right)=
\sum_{i=1}^{k}( \beta^{+}_{i,T}-(\beta')^{+}_{i,T} )=2k\varepsilon_{1}=(2^{\isB}k)(2^{1-\isB}\varepsilon_{1}).
\end{equation*}

To compute $S^{-}$, we know that  $\lambda'_{t}=\lambda_{t+1}$ for $1\leqslant t \leqslant r-1=r'$ (Proposition \eqref{prop:doublepartitiontypes}). Using Lemma \ref{lema:lagran_perm}, the position of $\overline{\lambda_{1}}=\overline{1}$ in $\wmax$ and the position of $\lambda_{1}=1$ in $\wprimemax$ are equal to $r$. In other words, $\wmax(r) = \overline{1}$, $\wprimemax(r) = 1$, and $\wmax(i)=\wprimemax(i)$ for every $i$ different from $r$.

For $i=r$, we have that $[i,i-1+\lambda_{r-i+1}]=[r,r]$. Thus, by Proposition \ref{prop:invprop_neg}, $\invn(\wmax) = \{(i,j) \tq  i\in [r] \mbox{ and } j\in [i , i-1+\lambda_{r-i+1}]\}$ and $\invn(\wprimemax) = \invn(\wmax)-\{(r,r)\}$. Assume that $(\widetilde\beta')^{-}_{r,r}=0$ to simplify our notation. 

For $(i,j)\in \invn(\wmax)$
%, i.e., $i \in [r]$ and $j\in [i, i-1+\lambda_{r-i+1}]$, all roots $\widetilde{\beta}_{i,j}^{-}$ and $(\widetilde{\beta}')_{i,j}^{-}$  are given as follows
\begin{align*}
\widetilde\beta^{-}_{i,j} &=
\left\{
\begin{array}{cl}
2^{1-\isB}\varepsilon_{1} & \mbox{if } i=r\\
-\varepsilon_{\wmax(i)}+\varepsilon_{1} & \mbox{if } i< r, j=r\\ 
2^{-\isB\cdot \delta_{ij}}(- \varepsilon_{\wmax(i)}-\varepsilon_{\wmax(j)}) & \mbox{otherwise}
\end{array} 
\right.\\
(\widetilde\beta')^{-}_{i,j} &= 
\left\{
\begin{array}{cl}
0 & \mbox{if } i=r\\
-\varepsilon_{\wmax(i)}-\varepsilon_{1} & \mbox{if } i< r, j=r\\ 
2^{-\isB\cdot \delta_{ij}}(- \varepsilon_{\wmax(i)}-\varepsilon_{\wmax(j)}) & \mbox{otherwise}
\end{array} 
\right.
\end{align*}
Then, $S^{-}$ can be rearranged as
\begin{equation*}
S^{-} = \sum_{i=1}^{r} \left( \sum_{j=i}^{i-1+\lambda_{r-i+1}} ( \widetilde\beta^{-}_{i,j}-(\widetilde\beta')^{-}_{i,j} ) \right) 
= \sum_{i=1} ^{r-1} (\widetilde\beta^{-}_{i,r}-(\widetilde\beta')^{-}_{i,r} ) + (\widetilde\beta^{-}_{r,r}-(\widetilde\beta')^{-}_{r,r} ) 
= (2^{\isB}r-2^{\isB}+1)(2^{1-\isB}\varepsilon_{1}).
\end{equation*}

We can easily observe that $2^{\isB}=\isB+1$. Therefore,
\begin{equation*}
\sigma(w)-\sigma(w') = (2^{\isB}k +2^{\isB}r-2^{\isB}+1)(2^{1-\isB}\varepsilon_{1}) = ((\isB +)T -\isB)(2^{1-\isB}\varepsilon_{1}). \qedhere
\end{equation*}
\end{proof}

The next two lemmas are required for the other proofs.  Consider $w\in \weyl^{(k)}_{n}$, with double partition $\Lambda=(\alpha,\lambda)$ and $r=\ell(\lambda)$.

\begin{lem}\label{lema:lambda_strict}
If $1\leqslant i<j\leqslant r$ then $j+\lambda_{r-j+1}\leqslant i+\lambda_{r-i+1}$.
\end{lem}
\begin{proof}
For every $i\in[r-1]$, we have $\lambda_{r-i+1}\geqslant 1+\lambda_{r-i}$ since $\lambda$ is a strict partition. Furthermore, if $l\in[r-i]$ then $\lambda_{r-i+1}\geqslant l+\lambda_{r-i+1-l}$. Hence, given $i<j$, take $l=j-i$.
\end{proof}

\begin{lem}\label{lema:conjugation}
Let $i$ and $j$ be integers such that $i \in [k]$ and $j\in [n-k]$. Then, $j\leqslant \alpha_{i}$ if, and only if, $k-i+1\leqslant \alpha_{j}^{*}$.
\end{lem}
\begin{proof} Notice that $j\leqslant \alpha_{i}$ if, and only if, $\{\alpha_{l}\colon \alpha_{l}\geq j\} \supseteq \{\alpha_{i},\ldots,\alpha_{k}\}$. The result follows from the definition of $\alpha^{*}$.
\end{proof}

\begin{prop}\label{prop:mainB2}
Let $w,w'$ be in $\weyl^{(k)}_{n}$ such that $\typew(w,w')=$ B2. Denote $T=w^{-1}(\overline{\lambda_{t}})=k+r-t+1$ and $Q=w^{-1}(v_{q})=k+r+q$. Then, $\kappa = T + Q -\isB$ and $\gamma = \varepsilon_{\lambda_{t}}-\varepsilon_{v_{q}}$.
\end{prop}
\begin{proof}
Observe that the indexes $T$ and $Q$ were chosen such that $v_{q}=\lambda_{t}-1$, $w(T) =\overline{\lambda_{t}}$, $
w(Q) =v_{q}$, $w'(T) =\overline{v_{q}}$, $w'(Q) = \lambda_{t}$, and $w(i)=w'(i)$ for every $i$ different from $T$ and $Q$. Recall that $\alpha_{i}'=\alpha_{i}$ for every $i$ (Proposition \eqref{prop:doublepartitiontypes}). Then, by Proposition \ref{prop:invprop_pos}, $\invp(w)= \invp(w')= \{(i,j) \tq i\in [k] \mbox{ and } j\in [k+1, k+\alpha_{i}]\}$.

Notice that $T \in [k+1,k+r]\subset [k+1,k+\alpha_{i}]$ and $Q\not\in [1,k+r]$. For $(i,j)\in\invp(w)$,
\begin{equation*}
\beta^{+}_{i,j}  = % \varepsilon_{w(i)}-\varepsilon_{w(j)} =
\left\{
\begin{array}{cl}
\varepsilon_{w(i)}+\varepsilon_{\lambda_{t}} & \mbox{if } j=T\\ 
\varepsilon_{w(i)}-\varepsilon_{v_{q}} & \mbox{if } j=Q\\ 
 \varepsilon_{w(i)}-\varepsilon_{w(j)} & \mbox{otherwise} %\mbox{if } j\neq T \mbox{ and } j\neq Q
\end{array} 
\right. ; \quad
(\beta')^{+}_{i,j}  = % \varepsilon_{w'(i)}-\varepsilon_{w'(j)} =
\left\{
\begin{array}{cl}
\varepsilon_{w(i)}+\varepsilon_{v_{q}} & \mbox{if } j=T\\ 
\varepsilon_{w(i)}-\varepsilon_{\lambda_{t}} & \mbox{if } j=Q\\ 
 \varepsilon_{w(i)}-\varepsilon_{w(j)} & \mbox{otherwise} %\mbox{, if } j\neq T \mbox{ and } j\neq Q
\end{array} 
\right.
\end{equation*}

The summation $S^{+}$ can be rearranged as follows
\begin{equation*}
S^{+}=  \sum_{i=1}^{k}\left( \sum_{j=k+1}^{k+\alpha_{i}}\left( \beta^{+}_{i,j}-(\beta')^{+}_{i,j} \right)\right).
\end{equation*}

Applying Lemma \ref{lema:conjugation}, $Q\leqslant k+\alpha_{i}$ if, and only if, $k-\alpha_{Q-k}^{*}+1 \leqslant i$.  Keeping in mind that $T$ always lies in the interval $[k+1,k+\alpha_{i}]$ for any $i$, we can split the summation $S^{+}$ over $i$ in the following parts: 
\begin{enumerate}
\item[(i)] If $i \in [k-\alpha_{Q-k}^{*}]$ then $Q\not\in[k+1,k+\alpha_{i}]$ and $\displaystyle\sum_{\mathclap{j=k+1}}^{k+\alpha_{i}} ( \beta^{+}_{i,j}-(\beta')^{+}_{i,j} ) = ( \beta^{+}_{i,T}-(\beta')^{+}_{i,T} ) = \varepsilon_{\lambda_{t}}-\varepsilon_{v_{q}}$.

\item[(ii)] If $i\in [k-\alpha_{Q-k}^{*}+1, k]$ then $Q\in[k+1,k+\alpha_{i}]$ and $\displaystyle\sum_{\mathclap{j=k+1}}^{k+\alpha_{i}} ( \beta^{+}_{i,j}-(\beta')^{+}_{i,j} ) = ( \beta^{+}_{i,T}-(\beta')^{+}_{i,T} )+ ( \beta^{+}_{i,Q}-(\beta')^{+}_{i,Q} ) = 2(\varepsilon_{\lambda_{t}}-\varepsilon_{v_{q}})$.
\end{enumerate}

Hence, $S^{+}$ is
\begin{equation*}
S^{+}= 
\sum_{i=1}^{k-\alpha_{Q-k}^{*}} \left( \sum_{j=k+1}^{k+\alpha_{i}} ( \beta^{+}_{i,j}-(\beta')^{+}_{i,j} ) \right)+
\sum_{i=k-\alpha_{Q-k}^{*}+1}^{k}\left( \sum_{ j=k+r+1}^{k+\alpha_{i}} ( \beta^{+}_{i,j}-(\beta')^{+}_{i,j} )\right)
= (k+\alpha_{Q-k}^{*})(\varepsilon_{\lambda_{t}}-\varepsilon_{v_{q}}).
\end{equation*}

To compute $S^{-}$, we know that $\lambda'_{t}=\lambda_{t}-1$ and $\lambda'_{i}=\lambda_{i}$ for $i\neq t$ (Proposition \eqref{prop:doublepartitiontypes}). Using Lemma \ref{lema:lagran_perm}, the position of $\overline{\lambda_{t}}$ in $\wmax$ is $\widetilde{T}:=\wmax^{-1}(\overline{\lambda_{t}})=T-k=r-t+1$, and the position of $v_{q}$ in $\wmax$ is $\widetilde{Q}:=\wmax^{-1}(v_{q})=Q-\alpha_{Q-k}^{*}$. In other words, $\wmax(\widetilde{T}) = w(T) = \overline{\lambda_{t}}$, $\wprimemax(\widetilde{T}) = w'(T) = \overline{v_{q}}$, $\wmax(\widetilde{Q}) = w(Q) = v_{q}$, $ \wprimemax(\widetilde{Q}) = w'(Q) = \lambda_{t}$,
and $\wmax(i)=\wprimemax(i)$ for every $i$ different from $\widetilde{T}$ and $\widetilde{Q}$.
By Proposition \ref{prop:invprop_neg}, $\invn(\wmax) = \{(i,j) \tq i\in [r] \mbox{ and } j\in [i , i-1+\lambda_{r-i+1}]\}$ and $\invn(\wprimemax) = \invn(\wmax)-\{(\widetilde{T},\widetilde{Q})\}$.

Notice that $\widetilde{T}\in [1,r]$ and $\widetilde{Q}\not\in [1,r]$. It remains to determine whether $\widetilde{T}$ and $\widetilde{Q}$ lie in the interval $[i,i-1+\lambda_{r-i+1}]$. Applying Lemma \ref{lema:conjugationExplict}, we have that $\widetilde{Q}=Q-\alpha_{Q-k}^{*}= k+r+q-\mu_{q}^{*}$. The number $\mu_{q}^{*}$ can be computed as it follows: $\mu_{q}^{*}= \#\{u_{l}\tq u_{l}>v_{q}\} = n-v_{q}-\#\{\lambda_{l} \tq \lambda_{l}>v_{q}\} - \#\{v_{l} \tq v_{l}>v_{q}\} = n - (\lambda_{t}-1)-(r-t+1)-(n-k-r-q)= k+t+q-\lambda_{t}$. Hence, $\widetilde{Q}=r-t+\lambda_{t}=\widetilde{T}-1+\lambda_{r-\widetilde{T}+1}$. 

Now, given $i\in[1,r]$, we can say whether $\widetilde{T}$ and $\widetilde{Q}$ belongs to $[i,i-1+\lambda_{r-i+1}]$ according to how $i$ compares to $\widetilde{T}$ by using Lemma \ref{lema:lambda_strict}. Namely,
\begin{align}
&\widetilde{T},\widetilde{Q}\in [i, i-1+\lambda_{r-i+1}], \mbox{ for } 1\leqslant i<\widetilde{T};\label{eq:B2_leqT}\\
&[\widetilde{T},\widetilde{Q}]= [i, i-1+\lambda_{r-i+1}], \mbox{ for } i=\widetilde{T};\label{eq:B2_eqT}\\
&\widetilde{T},\widetilde{Q}\not\in [i, i-1+\lambda_{r-i+1}], \mbox{ for } \widetilde{T}<i\leqslant r.\label{eq:B2_geqT}
\end{align}

Once $(\widetilde{T},\widetilde{Q})$ is not in $\invn(\wprimemax)$, we will assume that $(\widetilde\beta')^{-}_{\widetilde{T},\widetilde{Q}}=0$ to simplify our notation. For $(i,j)\in\invn(w)$,
\begin{align*}
\widetilde\beta^{-}_{i,j} & =
\left\{
\begin{array}{cl}
\varepsilon_{\lambda_{t}}-\varepsilon_{v_{q}} & \mbox{if } i=\widetilde{T}, j = \widetilde{Q}\\
2^{1-\isB}\varepsilon_{\lambda_{t}} & \mbox{if } i=\widetilde{T}, j=\widetilde{T}\\
\varepsilon_{\lambda_{t}}-\varepsilon_{\wmax(j)} & \mbox{if } i=\widetilde{T}, \widetilde{T}< j < \widetilde{Q}\\
- \varepsilon_{\wmax(i)}+\varepsilon_{\lambda_{t}} & \mbox{if } i<\widetilde{T}, j = \widetilde{T}\\
- \varepsilon_{\wmax(i)}-\varepsilon_{v_{q}} & \mbox{if } i<\widetilde{T}, j = \widetilde{Q}\\
2^{-\isB\cdot \delta_{ij}}(- \varepsilon_{\wmax(i)}-\varepsilon_{\wmax(j)}) & \mbox{otherwise}
\end{array} 
\right. \\
(\widetilde\beta')^{-}_{i,j}  & = 
\left\{
\begin{array}{cl}
0 & \mbox{, if } i=\widetilde{T}, j = \widetilde{Q}\\
2^{1-\isB}\varepsilon_{v_{q}} & \mbox{if } i=\widetilde{T}, j=\widetilde{T}\\
\varepsilon_{v_{q}}-\varepsilon_{\wmax(j)} & \mbox{if } i=\widetilde{T}, \widetilde{T}<j< \widetilde{Q}\\
- \varepsilon_{\wmax(i)}+\varepsilon_{v_{q}} & \mbox{if } i<\widetilde{T}, j = \widetilde{T}\\
- \varepsilon_{\wmax(i)}-\varepsilon_{\lambda_{t}} & \mbox{if } i<\widetilde{T}, j = \widetilde{Q}\\
2^{-\isB\cdot \delta_{ij}}(- \varepsilon_{\wmax(i)}-\varepsilon_{\wmax(j)}) & \mbox{otherwise}
\end{array} 
\right.
\end{align*}

The summation $S^{-}$ can be rearranged as follows
\begin{equation*}
S^{-} = \sum_{i=1}^{r} \left( \sum_{j=i}^{i-1+\lambda_{r-i+1}} \left( \widetilde\beta^{-}_{i,j}-(\widetilde\beta')^{-}_{i,j} \right) \right).
\end{equation*}

We can split the above summation over $i$ in the following parts:
\begin{enumerate}
\item[(i)] If $i \in [\widetilde{T}-1]$ then, by Equation \eqref{eq:B2_leqT}, 
\begin{equation*}
\sum_{j=i}^{\mathclap{i-1+\lambda_{r-i+1}}} ( \widetilde\beta^{-}_{i,j}-(\widetilde\beta')^{-}_{i,j} ) = ( \widetilde\beta^{-}_{i,\widetilde{T}}-(\widetilde\beta')^{-}_{i,\widetilde{T}} ) + ( \widetilde\beta^{-}_{i,\widetilde{Q}}-(\widetilde\beta')^{-}_{i,\widetilde{Q}} ) = 2(\varepsilon_{\lambda_{t}}-\varepsilon_{v_{q}}).
\end{equation*}
\item[(ii)] If $i=\widetilde{T}$ then, by Equation \eqref{eq:B2_eqT}, 
\begin{align*}
\sum_{\mathclap{j=\widetilde T}}^{\widetilde Q} ( \widetilde\beta^{-}_{\widetilde T,j}-(\widetilde\beta')^{-}_{\widetilde T,j}) &= 
( \widetilde\beta^{-}_{\widetilde T,\widetilde T}-  (\widetilde\beta')^{-}_{\widetilde T,\widetilde T})+  \sum_{\mathclap{j=\widetilde T +1}}^{\widetilde Q-1}( \widetilde\beta^{-}_{\widetilde{T},j}-(\widetilde\beta')^{-}_{\widetilde{T},j} )
 + ( \widetilde\beta^{-}_{\widetilde T,\widetilde Q}-(\widetilde\beta')^{-}_{\widetilde T,\widetilde Q})\\
& =(2^{1-\isB} +\widetilde{Q}-\widetilde{T})(\varepsilon_{\lambda_{t}}-\varepsilon_{v_{q}}).
\end{align*}

\item[(iii)] If $i\in [\widetilde{T}+1, r]$ then, by Equation \eqref{eq:B2_geqT}, $\displaystyle \sum_{j=i}^{{i-1+\lambda_{r-i+1}}} ( \widetilde\beta^{-}_{i,j}-(\widetilde\beta')^{-}_{i,j} ) = 0$.
\end{enumerate}

We can easily observe that $2^{1-\isB}-2=-\isB$. Thus, $S^{-}$ is
\begin{align*}
S^{-} & =  \sum_{i=1}^{\widetilde T-1}\left( \sum_{j=i}^{i-1+\lambda_{r-i+1}} ( \widetilde\beta^{-}_{i,j}-(\widetilde\beta')^{-}_{i,j} )\right) +
\sum_{\mathclap{j=\widetilde T}}^{\widetilde Q} ( \widetilde\beta^{-}_{\widetilde T,j}-(\widetilde\beta')^{-}_{\widetilde T,j})
+\sum_{\mathclap{i=\widetilde T+1}}^{r}\left( \sum_{j=i}^{i-1+\lambda_{r-i+1}} ( \widetilde\beta^{-}_{i,j}-(\widetilde\beta')^{-}_{i,j} )\right)  \\
&=  (-\isB+T+Q-k-\alpha_{Q-k}^{*})(\varepsilon_{\lambda_{t}}-\varepsilon_{v_{q}})
\end{align*}

Therefore,
\begin{equation*}
\sigma(w)-\sigma(w')  = (T+Q-\isB)(\varepsilon_{\lambda_{t}}-\varepsilon_{v_{q}}). \qedhere
\end{equation*}
\end{proof}

\begin{prop}\label{prop:mainB3}
Let $w,w'$ be in $\weyl^{(k)}_{n}$ such that $\typew(w,w')=$ B3. Denote $P=w^{-1}(u_{p})=p$ and $Q=w^{-1}(v_{q})=k+r+q$. Then, $\kappa = Q-P$ and $\gamma = \varepsilon_{u_{p}}-\varepsilon_{v_{q}}$.
\end{prop}
\begin{proof}
Observe that the indexes $P$ and $Q$ were chosen such that $v_{q}=u_{p}-x$, $w(P)=u_{p}$, $w(Q)=v_{q}$, $w'(P)=v_{q}$, $w'(Q)=u_{p}$, and $w(i)=w'(i)$ for every $i$ different from $P$ and $Q$. Recall that $\alpha'_{P}= \alpha_{P}-1$ and $\alpha'_{i}=\alpha_{i}$ for $i\neq P$ (Proposition \eqref{prop:doublepartitiontypes}). Then, by Proposition \ref{prop:invprop_pos}, $\invp(w)=\{(i,j) \tq 1\leqslant i\leqslant k \mbox{ and } k+1 \leqslant j \leqslant k+\alpha_{i}\}$ and $\invp(w')= \invp(w) - \{(P,Q)\}$.

Clearly $P\in[1,k]$ and $Q\not\in[1,k]$. Once $(P,Q)$ is not in $\invn(\wprimemax)$, we will assume that $(\widetilde\beta'_{P,Q})=0$. For $(i,j)\in\invp(w)$,
\begin{equation*}
\beta^{+}_{i,j}  = % \varepsilon_{w(i)}-\varepsilon_{w(j)} =
\left\{
\begin{array}{cl}
\varepsilon_{u_{p}} - \varepsilon_{v_{q}} & \mbox{if } i=P, j=Q\\ 
\varepsilon_{u_{p}}-\varepsilon_{w(j)} & \mbox{if } i= P, j\neq Q\\
\varepsilon_{w(i)}-\varepsilon_{v_{q}} & \mbox{if } i\neq P, j=Q\\ 
 \varepsilon_{w(i)}-\varepsilon_{w(j)} & \mbox{otherwise} %\mbox{, if } j\neq T \mbox{ and } j\neq Q
\end{array} 
\right.;
\quad
(\beta')^{+}_{i,j}  = % \varepsilon_{w'(i)}-\varepsilon_{w'(j)} =
\left\{
\begin{array}{cl}
0 & \mbox{, if } i=P, j=Q\\ 
\varepsilon_{v_{q}}-\varepsilon_{w(j)} & \mbox{if } i= P, j\neq Q\\
\varepsilon_{w(i)}-\varepsilon_{u_{p}} & \mbox{if } i\neq P, j=Q\\ 
 \varepsilon_{w(i)}-\varepsilon_{w(j)} & \mbox{otherwise} %\mbox{, if } j\neq T \mbox{ and } j\neq Q
\end{array} 
\right.
\end{equation*}

 Then, $S^{+}$ can be rearranged as follows
\begin{equation*}
S^{+}=  \sum_{i=1}^{k}\left( \sum_{j=k+1}^{k+\alpha_{i}} ( \beta^{+}_{i,j}-(\beta')^{+}_{i,j} )\right).
\end{equation*}

To compute this summation, we need to figure out when $Q\leqslant k+\alpha_{i}$. Notice that if $i=P$ then, by Equation \eqref{eq:alpha2}, $\alpha_{P}=r+\#\{v_{l}\tq v_{l}<u_{P}\}=r+q$, implying that $Q=k+\alpha_{P}$. Applying Lemma \ref{lema:conjugation}, $Q\leqslant k+\alpha_{i}$ if, and only if, $k-\alpha_{Q-k}^{*}+1 \leqslant i$.  But, $\alpha_{Q-k}^{*}=\alpha_{\alpha_{P}}^{*}=\#\{\alpha_{l}\tq \alpha_{l}>\alpha_{P}\}=k-P+1$. Hence, $Q\leqslant k+\alpha_{i}$ if, and only if, $i\geqslant P$.

We can split the summation of $S^{+}$ over $i$ in the following parts: 
\begin{enumerate}
\item[(i)] If $1\leqslant i < P$ then $Q\not\in [k+1, k+\alpha_{i}]$ and $\displaystyle\sum_{j=k+1}^{k+\alpha_{i}} ( \beta^{+}_{i,j}-(\beta')^{+}_{i,j} )=0$.

\item[(ii)] If $i = P$ then $[k+1, k+\alpha_{P}]=[k+1,Q]$ and $\displaystyle\sum_{\mathclap{j=k+1}}^{Q} ( \beta^{+}_{P,j}-(\beta')^{+}_{P,j} )= ( \beta^{+}_{P,Q}-(\beta')^{+}_{P,Q} )+\sum_{\mathclap{j=k+1}}^{Q-1} ( \beta^{+}_{P,j}-(\beta')^{+}_{P,j} )= (Q-k)(\varepsilon_{u_{p}}-\varepsilon_{v_{q}})$.

\item[(iii)] If $P< i \leqslant k$ then $Q\in [k+1, k+\alpha_{i}]$ and $\displaystyle\sum_{j=k+1}^{k+\alpha_{i}} ( \beta^{+}_{i,j}-(\beta')^{+}_{i,j} )= ( \beta^{+}_{i,Q}-(\beta')^{+}_{i,Q} )=\varepsilon_{u_{p}}-\varepsilon_{v_{q}}$.

\end{enumerate}

Hence, $S^{+}$ is
\begin{align*}
S^{+}&= \sum_{i=1}^{P-1} \left(\sum_{j=k+1}^{k+\alpha_{i}} ( \beta^{+}_{i,j}-(\beta')^{+}_{i,j} ) \right) +
\sum_{\mathclap{j=k+1}}^{Q} ( \beta^{+}_{P,j}-(\beta')^{+}_{P,j} )
 + \sum_{i=P+1}^{k} \left(\sum_{j=k+1}^{k+\alpha_{i}} ( \beta^{+}_{i,j}-(\beta')^{+}_{i,j} ) \right)\\
&= (Q-P)(\varepsilon_{u_{p}}-\varepsilon_{v_{q}}).
\end{align*}

To compute $S^{-}$, it is clear that both $\wmax$ and $\wprimemax$ are equal. Then,  the set of roots $\beta^{-}$ coincides with the set of roots $(\beta')^{-}$, and the sum $S^{-}$ should be zero. Therefore, 
\begin{equation*}
\sigma(w)-\sigma(w') = (Q-P)(\varepsilon_{u_{p}}-\varepsilon_{v_{q}}). \qedhere
\end{equation*}
\end{proof}

\begin{prop}\label{prop:mainB4}
Let $w,w'$ be in $\weyl^{(k)}_{n}$ such that $w,w'$ is a pair of type B4. Denote $P=w^{-1}(u_{p})=p$ and $T=w^{-1}(\overline{\lambda_{t}})=k+r-t+1$. Then, $\kappa = P+T-\isB$ and $\gamma = \varepsilon_{\lambda_{t}}-\varepsilon_{u_{p}}$.
\end{prop}
\begin{proof}
Observe that the indexes $P=p$ and $T=k+r-t+1$ were chosen such that $u_{p}=\lambda_{t}-x$, $w(P) =u_{p}$, $ w(T) = \overline{\lambda_{t}}$, $w'(P) =\lambda_{t}$, $ w'(T) =\overline{u_{p}}$, and $w(i)=w'(i)$ for every $i$ different from $P$ and $T$. Recall that $\alpha'_{P}= \alpha_{P}+x-1$ and $\alpha'_{i}=\alpha_{i}$ for $i\neq P$ (Proposition \eqref{prop:doublepartitiontypes}). Then, by Proposition \ref{prop:invprop_pos}, $\invp(w)=\{(i,j) \tq 1\leqslant i\leqslant k \mbox{ and } k+1 \leqslant j \leqslant k+\alpha_{i}\}$ and $ \invp(w') = \invp(w)\cup A$, where $A$ is the set given by: if $x>1$ then $A=\{(P,k+\alpha_{P}+l)\tq 1\leqslant l \leqslant x-1\}$; if $x=1$ then $A=\emptyset$.

Notice that $P\in[1,k]$ and $T\in [k+1, k+\alpha_{i}]$ for every $i\in [1,k]$ since $T \leqslant k+r\leqslant k+\alpha_{i}$. Then, all roots $\beta^{+}_{i,j}$ and $(\beta')^{+}_{i,j}$ for $(i,j)\in\invp(w)$, 
%i.e., $1\leqslant i \leqslant k$ and $k+1 \leqslant j \leqslant k+\alpha_{i}$, are given as follows
\begin{equation*}
\beta^{+}_{i,j}  = % \varepsilon_{w(i)}-\varepsilon_{w(j)} =
\left\{
\begin{array}{cl}
\varepsilon_{u_{p}} + \varepsilon_{\lambda_{t}} & \mbox{if } i=P, j=T\\ 
\varepsilon_{w(i)}+\varepsilon_{\lambda_{t}} & \mbox{if } i\neq P, j=T\\ 
\varepsilon_{u_{p}}-\varepsilon_{w(j)} & \mbox{if } i= P, j\neq T\\
 \varepsilon_{w(i)}-\varepsilon_{w(j)} & \mbox{otherwise}
\end{array} 
\right.
;\quad
(\beta')^{+}_{i,j}  = % \varepsilon_{w'(i)}-\varepsilon_{w'(j)} =
\left\{
\begin{array}{cl}
\varepsilon_{\lambda_{t}} + \varepsilon_{u_{p}} & \mbox{if } i=P, j=T\\ 
\varepsilon_{w(i)}+\varepsilon_{u_{p}} & \mbox{if } i\neq P, j=T\\ 
\varepsilon_{\lambda_{t}}-\varepsilon_{w(j)} & \mbox{if } i= P, j\neq T\\
 \varepsilon_{w(i)}-\varepsilon_{w(j)} & \mbox{otherwise}
\end{array} 
\right.
\end{equation*}

For $(P,j)\in A$, i.e., $k+\alpha_{P}+1\leqslant j \leqslant k+\alpha_{P}+x-1$,  the additional roots $(\beta')^{+}_{P,j}$ of $w'$ are $(\beta')^{+}_{P,j}= \varepsilon_{\lambda_{t}}-\varepsilon_{w(j)}$.

Then, $S^{+}$ can be rearranged as follows
\begin{equation*}
S^{+}=  \sum_{i=1}^{k}\left( \sum_{j=k+1}^{k+\alpha_{i}} ( \beta^{+}_{i,j}-(\beta')^{+}_{i,j} )\right) - \sum_{j=k+\alpha_{P}+1}^{k+\alpha_{P}+x-1} (\beta')^{+}_{P,j}.
\end{equation*}

Keeping in mind that $T$ always lies in the interval $[k+1,k+\alpha_{i}]$ for any $i$, we can split the above summation on $i$ in the following parts:
\begin{enumerate}
\item[(i)] If $1\leqslant i \leqslant k$ and $i\neq P$ then $\displaystyle\sum_{j=k+1}^{k+\alpha_{i}} ( \beta^{+}_{i,j}-(\beta')^{+}_{i,j} )=  \beta^{+}_{i,T}-(\beta')^{+}_{i,T} =\varepsilon_{\lambda_{t}}-\varepsilon_{u_{p}}$. 

\item[(ii)] If $i = P$ then $\displaystyle\sum_{\mathclap{j=k+1}}^{k+\alpha_{P}} ( \beta^{+}_{P,j}-(\beta')^{+}_{P,j} ) =  \sum_{\mathclap{j=k+1}}^{\mathclap{T-1}} ( \beta^{+}_{P,j}-(\beta')^{+}_{P,j} ) + ( \beta^{+}_{P,T}-(\beta')^{+}_{P,T} ) +\sum_{\mathclap{j=T+1}}^{k+\alpha_{P}} ( \beta^{+}_{P,j}-(\beta')^{+}_{P,j} ) =(1-\alpha_{P})(\varepsilon_{\lambda_{t}}-\varepsilon_{u_{p}})$.

\end{enumerate}

Hence, $S^{+}$ is
\begin{align*}
S^{+}&= \sum_{\substack{i \in[1,k] \\ i\neq T}} (\varepsilon_{\lambda_{t}}-\varepsilon_{u_{p}}) +(1-\alpha_{P})(\varepsilon_{\lambda_{t}}-\varepsilon_{u_{p}})
-\sum_{j=k+\alpha_{P}+1}^{k+\alpha_{P}+x-1} (\varepsilon_{\lambda_{t}}-\varepsilon_{w(j)})\\
& = (k-\alpha_{P})(\varepsilon_{\lambda_{t}}-\varepsilon_{u_{p}}) -\sum_{j=k+\alpha_{P}+1}^{k+\alpha_{P}+x-1} (\varepsilon_{\lambda_{t}}-\varepsilon_{w(j)}).
\end{align*}

To compute $S^{-}$, we know that $\lambda'_{t}=\lambda_{t}-x$ and $\lambda'_{i}=\lambda_{i}$ for $i\neq t$ (Proposition \eqref{prop:doublepartitiontypes}). Using Lemma \ref{lema:lagran_perm}, the position of $\overline{\lambda_{t}}$ in $\wmax$ is $\widetilde{T}:=\wmax^{-1}(w(T))=T-k=r-t+1$, and the position of $u_{p}$ in $\wmax$ is $\widetilde{P}:=\wmax^{-1}(w(P))=P+\alpha_{P}$. In other words, $\wmax(\widetilde{P}) = w(P)=u_{p}$, $\wprimemax(\widetilde{P}) = w'(P)=\lambda_{t}$, $
\wmax(\widetilde{T}) = w(T)=\overline{\lambda_{t}}$, $ \wprimemax(\widetilde{T}) = w'(T)=\overline{u_{p}}$, and $\wmax(i)=\wprimemax(i)$ for every $i$ different from $\widetilde{T}$ and $\widetilde{P}$. By Proposition \ref{prop:invprop_neg}, $\invn(\wmax) = \{(i,j) \tq 1\leqslant i\leqslant r \mbox{ and } i \leqslant j \leqslant i-1+\lambda_{r-i+1}\}$ and $\invn(\wprimemax) = \invn(\wmax)-B$, where $B=\{(\widetilde{T}, \widetilde{P} +l)\tq 0\leqslant l \leqslant x-1\}$.

Notice that $\widetilde{T}\in[1,r]$, $\widetilde{P}\not\in[1,r]$, and $\widetilde{P}=P+\alpha_{P}=u_{p}+\#\{\lambda_{l}\tq \lambda_{l}>u_{p}\} 
=(\lambda_{t}-x)+(r-t+1)=\widetilde{T}-x+\lambda_{r-\widetilde{T}+1}$.
Given $i\in[1,r]$, we can say whether $\widetilde{P}$ and $\widetilde{T}$ belongs to $[i, i-1+\lambda_{r-i+1}]$ according to how $i$ compares to $\widetilde{T}$ by applying Lemma \ref{lema:lambda_strict}. Namely,
\begin{align}
&\widetilde{T},\widetilde{P}\in [i, i-1+\lambda_{r-i+1}], \mbox{ for } 1\leqslant i<\widetilde{T},;\label{eq:B4_leqT}\\
&[\widetilde{T},\widetilde{P}+x-1]= [i, i-1+\lambda_{r-i+1}], \mbox{ for } i=\widetilde{T};\label{eq:B4_eqT}\\
&\widetilde{T},\widetilde{P}\not\in [i, i-1+\lambda_{r-i+1}], \mbox{ for } \widetilde{T}<i\leqslant r.\label{eq:B4_geqT}
\end{align}

Once $(\widetilde{T},j)\in B$, for $j\geqslant \widetilde{P}$, are not in $\invn(\wprimemax)$, we will assume that $(\widetilde\beta')^{-}_{\widetilde{T},j}=0$ to simplify our notation. 
For $(i,j)\in\invn(w)$,
\begin{align*}
\widetilde\beta^{-}_{i,j} &= 
\left\{ 
\begin{array}{cl}
\varepsilon_{\lambda_{t}}-\varepsilon_{u_{p}} & \mbox{if } i=\widetilde{T}, j = \widetilde{P}\\
2^{1-\isB}\varepsilon_{\lambda_{t}} & \mbox{if } i=\widetilde{T}, j=\widetilde{T}\\
\varepsilon_{\lambda_{t}}-\varepsilon_{\wmax(j)} & \mbox{if } i=\widetilde{T}, j> \widetilde{T}, j\neq \widetilde{P}\\
- \varepsilon_{\wmax(i)}+\varepsilon_{\lambda_{t}} & \mbox{if } i<\widetilde{T}, j = \widetilde{T}\\
- \varepsilon_{\wmax(i)}-\varepsilon_{u_{p}} & \mbox{if } i<\widetilde{T}, j = \widetilde{P}\\
 2^{-\isB\cdot \delta_{ij}}(- \varepsilon_{\wmax(i)}-\varepsilon_{\wmax(j)}) & \mbox{otherwise}
\end{array} 
\right. \\
(\widetilde\beta')^{-}_{i,j} &=
\left\{
\begin{array}{cl}
0 & \mbox{if } i=\widetilde{T}, j \geqslant \widetilde{P}\\
2^{1-\isB}\varepsilon_{u_{p}} & \mbox{if } i=\widetilde{T}, j=\widetilde{T}\\
\varepsilon_{u_{p}}-\varepsilon_{\wmax(j)} & \mbox{if } i=\widetilde{T}, \widetilde{T}< j < \widetilde{P}\\
- \varepsilon_{\wmax(i)}+\varepsilon_{u_{p}} & \mbox{if } i<\widetilde{T}, j = \widetilde{T}\\
- \varepsilon_{\wmax(i)}-\varepsilon_{\lambda_{t}} & \mbox{if } i<\widetilde{T}, j = \widetilde{P}\\
2^{-\isB\cdot \delta_{ij}}(- \varepsilon_{\wmax(i)}-\varepsilon_{\wmax(j)}) & \mbox{otherwise}
\end{array} 
\right.
\end{align*}

Then, $S^{-}$ can be rearranged as follows
\begin{equation*}
S^{-} = \sum_{i=1}^{r} \left( \sum_{j=i}^{i-1+\lambda_{r-i+1}} ( \widetilde\beta^{-}_{i,j}-(\widetilde\beta')^{-}_{i,j} ) \right).
\end{equation*}

We can split the above summation on $i$ in the following parts:
\begin{enumerate}
\item[(i)] If $1\leqslant i < \widetilde{T}$ then, by Equation \eqref{eq:B4_leqT}, $\displaystyle\sum_{j=i}^{\mathclap{i-1+\lambda_{r-i+1}}} ( \widetilde\beta^{-}_{i,j}-(\widetilde\beta')^{-}_{i,j} )= ( \widetilde\beta^{-}_{i,\widetilde{T}}-(\widetilde\beta')^{-}_{i,\widetilde{T}} )+( \widetilde\beta^{-}_{i,\widetilde{P}}-(\widetilde\beta')^{-}_{i,\widetilde{P}} )=2(\varepsilon_{\lambda_{t}}-\varepsilon_{u_{p}})$.

\item[(ii)] If $i = \widetilde{T}$ then, by Equation \eqref{eq:B4_eqT}, $\displaystyle\sum_{j=\widetilde{T}}^{\mathclap{\widetilde{P} +x-1}} ( \widetilde\beta^{-}_{\widetilde{T},j}-(\widetilde\beta')^{-}_{\widetilde{T},j} )= ( \widetilde\beta^{-}_{\widetilde{T},\widetilde{T}}-(\widetilde\beta')^{-}_{\widetilde{T},\widetilde{T}} ) + \sum_{j=\widetilde{T}+1}^{\mathclap{\widetilde{P}-1}} ( \widetilde\beta^{-}_{\widetilde{T},j}-(\widetilde\beta')^{-}_{\widetilde{T},j} ) + ( \widetilde\beta^{-}_{\widetilde{T},\widetilde{P}}-(\widetilde\beta')^{-}_{\widetilde{T},\widetilde{P}} ) 
+ \sum_{j=\widetilde{P}+1}^{\mathclap{\widetilde{P} +x-1}} ( \widetilde\beta^{-}_{\widetilde{T},j}-(\widetilde\beta')^{-}_{\widetilde{T},j} )
=  (2^{1-\isB} + \widetilde{P}-\widetilde{T})(\varepsilon_{\lambda_{t}}-\varepsilon_{u_{p}}) + \sum_{j=P+\alpha_{P}+1}^{P+\alpha_{P}+x-1}(\varepsilon_{\lambda_{t}}- \varepsilon_{\widetilde{w}(j)})$.

\item[(iii)] If $\widetilde{T}< i\leqslant r$ then, by Equation \eqref{eq:B4_geqT}, $\displaystyle\sum_{j=i}^{\mathclap{i-1+\lambda_{r-i+1}}} ( \widetilde\beta^{-}_{i,j}-(\widetilde\beta')^{-}_{i,j} )=0$.

\end{enumerate}

Hence,
\begin{equation*}
S^{-} = \sum_{i=1}^{\widetilde T-1}\left(\sum_{j=i}^{ i-1+\lambda_{r-i+1}}  ( \widetilde\beta^{-}_{i,j}-(\widetilde\beta')^{-}_{i,j} ) \right) +
\sum_{\mathclap{j=\widetilde T}}^{\mathclap{\widetilde P +x-1}} ( \widetilde\beta^{-}_{\widetilde T,j}-(\widetilde\beta')^{-}_{\widetilde T,j}) +
\sum_{\mathclap{i=\widetilde T+1}}^{r}\ \left( \sum_{j=i}^{ i-1+\lambda_{r-i+1}} ( \widetilde\beta^{-}_{i,j}-(\widetilde\beta')^{-}_{i,j} )\right).
\end{equation*}
Thus,
\begin{equation*}
S^{-}=(P+T+\alpha_{p}-k-\isB)(\varepsilon_{\lambda_{t}}-\varepsilon_{u_{p}}) +\sum_{j=P+\alpha_{P}+1}^{P+\alpha_{P}+x-1}(\varepsilon_{\lambda_{t}}- \varepsilon_{\wmax(j)}).
\end{equation*}

Finally, the remaining summation on $S^{-}$ is supposed to cancel with the one in $S^{+}$. If we prove that $\wmax(j)=w(j+k-P)$ for $j\in [P+\alpha_{P}+1,P+\alpha_{P}+x-1]$ then we clearly have
\begin{equation*}
\sum_{j=P+\alpha_{P+1}}^{P+\alpha_{P}+x-1}(\varepsilon_{\lambda_{t}}- \varepsilon_{\widetilde{w}(j)})= \sum_{j=P+\alpha_{P+1}}^{P+\alpha_{P}+x-1}(\varepsilon_{\lambda_{t}}- \varepsilon_{w(j+k-P)})= \sum_{j=k+\alpha_{P+1}}^{k+\alpha_{P}+x-1}(\varepsilon_{\lambda_{t}}- \varepsilon_{w(j)}).
\end{equation*}

Considering $l=j-P-\alpha_{P}$, the above assertion is equivalent to prove that $w(l+\alpha_{P}+k)=\wmax(l+P+\alpha_{P})$, for $l\in[1,x-1]$. Clearly, $l+\alpha_{P}+k> k$ and, by Lemma \ref{lema:lagran_perm}, $\wmax^{-1}(w(l+\alpha_{P}+k)) = l+\alpha_{P}+k - \alpha_{l+\alpha_{P}}^{*}$.
We know that $\alpha_{P+1}=\alpha'_{P+1}\geqslant \alpha_{P}' = \alpha_{P}+x-1$ since $\alpha'$ is also a partition. Then, by definition, $\alpha_{l+\alpha_{P}}^{*}=\#\{i \tq \alpha_{i}\geqslant l+\alpha_{P}\}=\#\{P+1, \dots, k\}=k-P$ for every $l\in [1,x-1]$. Hence, $\wmax^{-1}(w(l+\alpha_{P}+k))=l+P+\alpha_{P}$.

Therefore,
\begin{equation*}
\sigma(w)-\sigma(w') = (P+T-\isB)(\varepsilon_{\lambda_{t}}-\varepsilon_{v_{q}}). \qedhere
\end{equation*}
\end{proof}

\section{Final comments and further directions}\label{sec:finalcomments}

We conclude mentioning related facts and future research based on the results of this work.

\begin{enumerate}
\item We have figured out that Theorem \ref{thm:dualheight} is true for all flag manifolds of split real forms. This will be presented in a forthcoming paper.
\item Example \ref{ex:incidencegraph} of $\mathrm{OG}(2,9)$ had revealed how the choice of signals may be an obstacle to get the whole homology. We hope to apply the theory developed in \cite{RS19} and find an algorithm that provides a good choice for general covering pairs. 
\item For type D, if one obtains a complete characterization of the covering pairs, the same techniques may be applied to achieve similar results with the visual interpretations for the coefficients of even Orthogonal Grassmannians.
\item As well as we had established the HSYD's from the shifted Young diagrams defined by \cite{IN08} and \cite{GK15}, we expect to define analogous excitations inside the HSYD's for the general context of any isotropic Grassmannians.
\end{enumerate}

\section*{Acknowledgment}
We thank to David Anderson for helpful suggestions and valuable discussions on an earlier version.

\bibliographystyle{amsplain}
\bibliography{biblio}
\end{document}